\documentclass[oneside]{amsart}
\usepackage{preamble}
\title{The period map from commutative to noncommutative deformations}
\author{Samuel A.\ Moore}
\address{St. John's College, University of Oxford}
\email{samuel.moore@maths.ox.ac.uk}
\newcommand{\citelink}[2]{\hyperlink{cite.\therefsection @#1}{#2}}

\begin{document}

\begin{abstract}
We study the period map from infinitesimal deformations of a qcqs derived scheme $X$ over a field $k$ to those of the associated $k$-linear $\infty$-category $\QC(X)$. We identify the corresponding map on tangent fibres with the dual HKR map $\mrm{R}\Gamma(X, \mathbf{T}_{X/k})[1] \longrightarrow \HH^\bullet(X/k)[2]$ and study its injectivity on homotopy groups. As applications, we show liftability of qcqs schemes along square-zero extensions to be a derived invariant, at least when $\mrm{char}(k) \ne 2$, and exhibit cases where the entire classical deformation functor of $X$ is a derived invariant; this partially answers a question of Lieblich.
\end{abstract} 

\maketitle

\tableofcontents

\section{Introduction}

Let $X$ be a scheme over a field $k$ and $\Lambda$ be a complete Noetherian local ring with residue field $k$. The formal deformation theory of $X$ over Artinian local $\Lambda$-algebras with residue field $k$ is a classical object of study, but one can also study `noncommutative' variants, such as the deformation theory of the stable $k$-linear $\infty$-category $\QC(X)$ of quasicoherent sheaves; this has been developed by many authors, with early contributions including \cite{hom-mirror}, \cite{seidel}, \cite{lowen-vdb1} and \cite{kkp}.

This paper studies the relationship between these two deformation problems in the generality of derived algebraic geometry, where one may consider deformations of $X$ and $\QC(X)$ over animated (equivalently simplicial commutative) rings. There are corresponding derived deformation functors $\mrm{Def}_{X, \Lambda}$ and $\mrm{CatDef}_{\QC(X), \Lambda}$, which constitute \emph{formal moduli problems} in the sense of \cite{DAGX}, at least after suitably completing the latter functor. It is well-known that many infinitesimal deformation problems are controlled, at least in characteristic $0$, by differential graded Lie algebras (DGLA), cf.\ \cite{deligne-letter}, \cite{drinfeld-letter}, \cite{hinich-formal-stacks}. A theorem of Lurie and Pridham in equal characteristic $0$, and Brantner--Mathew more generally, states more precisely that the $\infty$-category of formal moduli problems over Artinian local animated $\Lambda$-algebras with residue field $k$ is equivalent to that of derived $(\Lambda,k)$-partition Lie algebras $\mrm{Alg}_{\mrm{Lie}^{\pi}_{\Lambda, \Delta}}$ (which are equivalently shifted DGLA when $\Lambda = k$ is a characteristic $0$ field) via the functor sending a formal moduli problem $F$ to its \emph{tangent fibre} $\mrm{T}_F$ \cite{pridham}, \cite{DAGX}, \cite{BM}.

In this paper, we construct and study the `commutative-to-noncommutative' period map \[\theta \colon \mrm{Def}_{X, \Lambda} \longrightarrow \mrm{CatDef}_{\QC(X), \Lambda},\] defined for derived schemes $X$ over $k$, which sends a lift $Y$ to an animated ring $A$ to the $A$-linear $\infty$-category $\QC(Y)$ lifting $\QC(X)$. On tangent fibres, this period map corresponds to a morphism 
\[\mrm{T}_\theta \colon \mrm{R}\Gamma(X, \mathbf{T}_{X/k})[1] \longrightarrow \mrm{HH}^\bullet(X/k)[2]\] of derived partition Lie algebras from the shifted derived global sections of the tangent complex to the twice-shifted Hochschild cochains.

Variants of such a period map appear in the literature, at least in characteristic $0$. For example, To\"{e}n and To\"{e}n--Anel mention a variant of $\theta$ for smooth and proper schemes in \cite[\S 3.3.\ (6)]{toen-higher} and \cite[\S 8(3)]{toenAnel}, and state its effect on tangent spaces (that is, they describe $\pi_0 \mrm{T}_{\theta}$), and the underlying map of DGLAs in characteristic 0 appears in \cite[Theorem 4.18]{kkp}.

However, a significant advantage of approaching the map $\theta$ via derived algebraic geometry and the formalism of formal moduli problems is that the latter possess a canonical, functorial obstruction theory. The map $\mrm{T}_{\theta}$ therefore automatically respects obstruction classes, which is crucial for the applications we have in mind. Moreover, $\mrm{T}_{\theta}$ still plays the same role in positive or mixed characteristic (where bare DGLA would not suffice) or for more general derived schemes.

\subsection{Statement of results}
Let us first rapidly recall some generalities on Hochschild homology for a qcqs derived scheme $X/k$. The Hochschild homology sheaf $\underline{\HH}_\bullet$ can be equipped with the HKR filtration $\mrm{Fil}^\star_{\mrm{HKR}} \underline{\HH}_{\bullet}$; when $X$ is smooth over $k$, this is simply the Whitehead filtration $\tau_{\ge \star} \underline{\HH}_\bullet$.  Should this filtration split, we say that the \emph{strong HKR theorem} holds for $X$. More generally, there is an associated spectral sequence, called the \emph{HKR spectral sequence}, with signature \[\mrm{E}_2^{i,j}=\mrm{H}^{j}(X, \Lambda^i \mrm{L}_{X/k}) \implies \HH_{i-j}(X/k),\] where the differential $d_r$ has bidegree $(r-1,r)$, and one says that the \emph{weak} HKR theorem holds for $X$ if this spectral sequence collapses on the $\mrm{E}_2$ page (this terminology is borrowed from \cite{AV}).

In any characteristic, there is an `HKR map' $\underline{\mrm{HH}}_\bullet \to \mrm{L}_{X/k}[1]$ which induces an equivalence $\mrm{gr}^1_{\mrm{HKR}} \underline{\mrm{HH}}_\bullet \xrightarrow{\simeq} \mrm{L}_{X/k}[1]$ (see \cref{relating}). We call the dual map \[\mrm{HKR}^\vee \colon \mrm{R}\Gamma(X, \mathbf{T}_{X/k})[-1] \longrightarrow \mrm{HH}^\bullet(X/k),\] landing in the Hochschild cochains, the \emph{dual} HKR map.

We may now state the main result of this paper: 
\begin{mainthm}\label{thm-A} Let $X$ be a qcqs derived $k$-scheme. The commutative-to-noncommutative period map \[\theta\colon \mrm{Def}_{X, \Lambda} \longrightarrow \mrm{CatDef}_{\QC(X), \Lambda}\] enjoys the following properties:
\begin{enumerate}
\item On tangent fibres, the induced map $\mrm{T}_{\theta}\colon \mrm{R}\Gamma(X, \mathbf{T}_{X/k})[1] \longrightarrow \HH^\bullet(X/k)[2]$ is the dual {HKR} map.
\item $\mrm{T}_{\theta}$ is injective on $\pi_i$:
\begin{enumerate} \item For \emph{all} $i$ if $\mrm{char}(k)=0$, and for $i \ge n+2-p$ if $\mrm{char}(k)=p > 0$ and $X$ is $n$-truncated.
\item \label{strong-I} For \emph{all} $i$ if the strong HKR theorem holds for $X$. \end{enumerate}
 \item \label{CY-I} Suppose moreover that $X$ is a weak Calabi--Yau variety of dimension $d$. Then $\mrm{T}_{\theta}$ is injective on $\pi_{i}$ if (and only if) all differentials in the HKR spectral sequence out of the terms $\{\mrm{E}^{1,d+i-1}_{r}\}_{r \ge 2}$ vanish. 
\end{enumerate}
\end{mainthm}

\begin{remark} The conditions on the HKR filtration appearing in \hyperref[strong-I]{(\ref*{strong-I})} and \hyperref[CY-I]{(\ref*{CY-I})} are often satisfied: the strong HKR theorem holds in characteristic $0$ or when $X$ is smooth with $\dim X < \mrm{char}(k)$ (see \cref{degen}), or if $X/k$ is smooth and proper with $\dim X \le p = \mrm{char}(k)$ \cite[Corollary 1.5.]{AV}.
\end{remark}

Supplementing the aforementioned results on the HKR spectral sequence, we also observe the following result of independent interest:
\begin{restatable}{proposition}{hkrcriterion} 	
\label{all-four-condition} 
Let $X$ be a smooth and proper scheme over a perfect field $k$ of characteristic $p > 0$, and suppose that 
\begin{enumerate}
	\item the crystalline cohomology $\mrm{H}^*_{\mrm{cris}}(X/W)$ is torsion-free, and
	\item the Hodge--de Rham spectral sequence for $X$ degenerates at $\mrm{E}_1$. 
\end{enumerate} 
Then the $\mrm{HKR}$, Tate, and de Rham--$\mrm{HP}$ spectral sequences for $X$ degenerate.
\end{restatable}

Our results have various ramifications for deformation theory, owing to the functoriality of obstruction classes attached to a formal moduli problem. We summarise this in the following: 
\begin{restatable}{mainthm}{thmB}\label{thm-B}
Let $X$ and $Y$ be qcqs derived $k$-schemes. Assume also that $\pi_{-1}\mrm{T}_{\theta}\colon \mrm{H}^2(\mrm{T}_{Z/k}) \to \mrm{HH}^3(Z/k)$ is injective for $Z = X$ and $Y$ (e.g. if $X,Y$ are classical and $\mrm{char}(k) \ne 2$).
\begin{enumerate} \item \label{minus-one} Suppose that $X$ has a lift $\mathfrak X/B$ for some local Artinian ring $B$ with residue field $k$. Then for any square-zero extension $A \to B$, $\mathfrak X$ lifts to $A$ if and only if $\QC(\mathfrak X)$ does. 
	
Hence if $\mathfrak X$ and $\mathfrak Y$ are derived-equivalent lifts of $X$ and $Y$ over $B$, then $\mathfrak X$ lifts to $A$ if and only if $\mathfrak Y$ does.
\item \label{abelianity} If $\mrm{T}_{\theta}$ is injective on $\pi_i$ for \emph{all} $i$, and the partition Lie algebra $\mrm{HH}^\bullet(X/k)[2]$ is abelian, then so is the partition Lie algebra $\mrm{R}\Gamma(X, \mathbf {T}_{X/k})[1]$. 
\end{enumerate}
\end{restatable}

Here we call two $R$-schemes $Z$ and $Z'$ \emph{derived equivalent} if there exists an $R$-linear equivalence of $\infty$-categories $\QC(Z) \simeq \QC(Z')$ (cf.\ \cref{derEquiv}). We in fact prove a more general result concerning square-zero extensions of animated rings; see \cref{lifts}. 

\begin{remark} The conclusion of \hyperref[minus-one]{Theorem \ref*{thm-B} (\ref*{minus-one})} does not hold for stacks in positive characteristic in general. Indeed, the stack $B \alpha_p$ over a perfect field $k$ of characteristic $p$ does not lift to $W_2(k)$, whilst the $\infty$-category $\QC(B \alpha_p) \simeq \Mod_{k[t]/t^p}$ has the evident lift $\Mod_{W(k)[t]/t^p}$ to $W(k)$.
\end{remark}

We remark that in positive characteristic, \hyperref[thm-B]{Theorem \ref*{thm-B} (1)} shows that the hypotheses of the Hodge--de Rham degeneration theorem of Deligne--Illusie \cite[Corollaire 2.4.]{DI} are derived invariants among smooth and proper $k$-schemes; this sharpens a theorem of Antieau--Bragg \cite[Theorem 2.6]{AntieauBragg} and clarifies some of the behaviour of the special derived-equivalent Calabi--Yau threefolds constructed by Addington--Bragg \cite{Add-Bragg}. 

Finally, we observe that in some special cases of interest, $\theta$ even induces an isomorphism of (classical) deformation functors:  
\begin{corollary}\label{def-iso} \,
\begin{enumerate} \item \label{cor1} Suppose $Z$ is a smooth and proper $k$-scheme, and suppose that the groups $\mrm{H}^1(Z, \mathcal O_Z)$, $\mrm{H}^2(Z, \mathcal O_Z)$, and $\mrm{H}^{0}(Z, \land^2\mrm{T}_{Z/k})$ vanish. Then the period map induces an isomorphism $\mrm{Def}_{Z, \Lambda}^{\mrm{cl}} \xrightarrow{\cong} \mrm{CatDef}^{\land, \mrm{cl}}_{\QC(Z), \Lambda}$ of classical deformation functors. 
\item \label{cor2} In particular, if $X$ and $Y$ are two such $k$-schemes, then any derived equivalence between $X$ and $Y$ gives rise to an isomorphism of classical deformation functors \[\mathrm{Def}_{X, \Lambda}^{\mathrm{cl}} \cong {\mathrm{Def}_{Y, \Lambda}^{\mathrm{cl}}}.\]
If $X$ and $Y$ are moreover projective, this induces a bijection between lifts of $X$ and $Y$ to $\Lambda$.
\end{enumerate}
\end{corollary}
\cref{def-iso} applies in particular to strict Calabi--Yau varieties of dimension at least three satisfying Hodge symmetry. Applied to such Calabi--Yau varieties in characteristic $0$, \cref{cor1} is essentially the same as \cite[Theorem 1.1]{morimura}, whilst \hyperref[cor2]{(\ref*{cor2})} was proven by Ward for strict Calabi--Yau threefolds in characteristic $p > 3$ satisfying stronger conditions \cite[Proposition 4.5.15]{Ward}. One novelty of our result is that no liftability is needed.

\subsection{Relation to other work} We note that recently, in independent unpublished work, Kurama--Perry have obtained many of the same results, including \hyperref[thm-B]{Theorem \ref*{thm-B}} (1), and \cref{def-iso} for smooth schemes.

Lieblich has asked how far liftability is a derived invariant of schemes in positive/mixed characteristic; see for example \cite[\S 5.1]{Ward}, \cite[Question 2.7.]{AntieauBragg}, \cite[Question 5.3.2]{Voet}, and \cite[\S 5.4]{Grigg} for related discussion. \hyperref[thm-B]{Theorem \ref*{thm-B}} and \cref{def-iso} (2) give positive answers in various cases. 

Hacon--Witaszek have shown (assuming some facts about resolutions of singularities) that if $k$ is a perfect field of characteristic $p > 5$, then liftability (to a given mixed characteristic complete DVR $R$ of residue field $k$) is a birational invariant within (smooth, projective) strict Calabi--Yau threefolds over $k$ \cite[Corollary 6.2]{Hacon-Witaszek}. In work to appear, Kurama--Perry have extended this result to $p > 3$ unconditionally by generalising Bridgeland's interpretation of flops via derived categories (see \cite{Bridgeland}) to mixed characteristic, and lifting flops to characteristic $0$ \cite{KuramaPerry}. Since birational Calabi--Yau threefolds are derived equivalent (at least under these assumptions), there is substantial overlap between these results and \cref{def-iso} (2); we are curious about their interaction.

Finally, let us mention a connection to Calabi--Yau varieties. A strong form of the Bogomolov--Tian--Todorov (BTT) theorem of \cite{Bog}, \cite{Tian}, and \cite{Todorov} states that the DGLA controlling deformations of a weak Calabi--Yau variety $X$ over a characteristic $0$ field $k$ is (homotopy) abelian; see \cite{Goldman-Millson} over $\mathbf C$, and \cite[Theorem (B)]{Iacono-Manetti} more generally. The analogous result concerning deformations of a smooth and proper Calabi--Yau $k$-linear $\infty$-category $\mathscr C$ holds in characteristic $0$, as sketched by Kontsevich--Katzarkov--Pandit \cite[Theorem 4.30]{kkp} (see also \cite[Theorem 2]{Terilla}, \cite{iwanari16}, \cite{iwanari-BTT}).  Upcoming work of Brantner--Devalapurkar--Horel proves a version of the noncommutative BTT theorem for some Calabi--Yau $\infty$-categories in positive characteristic.
\subsection{Outline of proof} 
Let us briefly outline the proof of \hyperref[thm-A]{Theorem \ref*{thm-A}}. If $\mrm{sqz}_k(V)$ denotes the square-zero extension of $k$ by some $V \in \Mod_{k, \ge 0}$, then it is known that
\[\Omega \mrm{Def}_{X, \Lambda}(\mrm{sqz}_k(V)) \simeq \mrm{Map}_{\QC(X)}(\mrm{L}_{X/k}, \mathcal O_X \otimes_k V),\] where we have taken the basepoint at the trivial deformation. We will construct similar equivalences
\[\Omega \mrm{CatDef}^{\land}_{\QC(X), \Lambda}(\mrm{sqz}_k(V)) \simeq \mrm{Map}_{\Fun_k^{\mrm{L}}(\QC(X), \QC(X))}(\mrm{id}, \mrm{id} \otimes_k V[1]) \simeq \mrm{Map}_{\QC(X)}(\underline{\mrm{HH}}_\bullet[-1], \mathcal O_X \otimes_k V),\] where the first map is a variant of the \emph{Atiyah class} construction.

We can view the right sides of these equivalences as functors on $\QC(X)_{\ge 0}$ evaluated at sheaves of the form $\mathcal O_X \otimes_k V$. We show that the `looped' period map accordingly extends to a transformation \[h \colon \mrm{Map}_{\QC(X)}(\mrm{L}_{X/k}, -) \longrightarrow \mrm{Map}_{\QC(X)}(\underline{\mrm{HH}}_\bullet[-1], -)\] of functors on $\QC(X)_{\ge 0}$. By generalising a result of Căldăraru, relating the universal Atiyah class to the HKR map, we identify the map $\underline{\mrm{HH}}_\bullet \to \mrm{L}_{X/k}[1]$ classifying $h$. The other results then follow by further analysis of the dual HKR map and results on formal moduli problems.

\subsection{Organisation}
\cref{defs} is mostly preparatory in nature: we review the construction of the two formal moduli problems involved, providing the necessary alterations which allow us to consider mixed characteristic deformations. \cref{hochschild-atiyah} studies the Atiyah class and its relation to categorical deformations. In \cref{period}, we construct the period map $\theta$, identify its effect on tangent fibres using the results of \S\S 2--3, and discuss conditions under which $\mrm{T}_{\theta}$ is injective on homotopy groups. In the final section, we explain more concretely what this entails for deformation theory, and develop some basic theory of \emph{abelian} partition Lie algebras.

\subsection{Notation and conventions}
We fix throughout a (discrete) complete Noetherian local ring $\Lambda$ with residue field $k$.

We employ $\infty$-categories throughout, primarily following Lurie's conventions in \cite{HTT}, \cite{HA}, and \cite{SAG} and adopting an `implicitly derived' convention; in particular $\mathcal S$ denotes the $\infty$-category of spaces (or $\infty$-groupoids, etc.) and $\Mod_{R}$ denotes the $\infty$-category of $R$-module spectra. We primarily use homological indexing, with the convention that $\mrm{H}^i = \pi_{-i} = \mrm{H}_{-i}$; accordingly, we also write $\underline{\pi}_i = \mathscr{H}^{-i}$ for the $i$th homotopy sheaf. We write $\mrm{CAlg}^{\mrm{an}}$ for the $\infty$-category of \emph{animated rings} (which may be modelled by simplicial commutative rings) and $\mrm{DSch}$ for the $\infty$-category of derived schemes in the sense of \cite[Definition 4.2.8]{DAGV} (as characterised by [\citelink{DAGV}{Ibid}, Theorem 4.2.15]). We write $\QC(Y)$ for the $\infty$-category of quasicoherent sheaves on a derived scheme $Y$ as defined (for the underlying spectral scheme) in \cite[Definition 2.2.2.1]{SAG}; if $Y$ is classical, its homotopy category identifies with $\mrm{D}_{\mrm{qc}}(Y)$ [\citelink{SAG}{Ibid}, \S 2.2.6.2]. Given a presentable stable $\infty$-category $\mathscr C$, we write $\mrm{map}_{\mathscr C}(-,-)$ for the mapping spectra between objects. Given a map $f \colon x \to y$ in an $\infty$-category $\mathscr C$, we write $\mathscr C_{x \sslash y}$ for $(\mathscr C_{/y})_{x/} \simeq (\mathscr C_{x/})_{/y}$. If $\mathscr C$ is an $\infty$-category, we denote its full subcategory of compact objects by $\mathscr C^\omega$.

By a weak Calabi--Yau variety over $k$, we mean an equidimensional smooth and proper $k$-scheme $X$ with trivial canonical bundle; a strict Calabi--Yau variety additionally has vanishing $\mrm{H}^i(X, \mathcal O_X)$ for $0 < i < \dim X$.

\subsection{Acknowledgements}
The author is very grateful to Lukas Brantner for his guidance and for suggesting this problem. We also thank Riku Kurama, Joshua Mundinger, Alexander Perry, Alexander Petrov and Joseph Stahl for helpful discussions.

\subsection{AI usage declaration} Some versions of ChatGPT Edu (up to 5.5) were used to find references and proofread, but not to generate text or proofs.

\section{Deformations and Hochschild (co)homology}\label{defs}

We begin by recalling some basic theory of formal moduli problems and constructing the examples relevant to the rest of the paper. As much of the literature focuses on equicharacteristic deformations (i.e.\ over augmented Artinian local $k$-algebras), we supply some modifications to accommodate mixed characteristic deformations.

\subsection{Generalities on formal moduli problems} \label{mixeddefs}\begin{notation}\label{a-scr} Throughout this subsection, we let $\mathscr A$ denote one of the $\infty$-categories $\CAlg^{\an}_{\Lambda \sslash k}, \CAlg^{\an}_{/k}$, $\mrm{Alg}^{\mathbf E_n}_{\Lambda \sslash k} = \mrm{Alg}^{\mathbf E_n }(\Mod_{\Lambda})_{/k}$, or $\mrm{Alg}^{\mathbf E_n}_{/k}$, where $0 \le n \le \infty$. In particular, $\mrm{Alg}^{\mathbf E_0}_{\Lambda \sslash k} = (\Mod_{\Lambda})_{\Lambda \sslash k}$.
\end{notation}

The trivial square zero extensions of $k$ by $k[m]$ assemble to a spectrum object $\{ \mrm{sqz}_k(k[m])\}_{m \in \N}$ of $\CAlg^{\an}_{k \sslash k}$. This forgets to a spectrum object in each of the other $\mathscr A$ considered above, and we may regard each $(\mathscr A, \{\mrm{sqz}_k(k[m])\}_{m \in \N})$ as a deformation context in the sense of \cite[Definition 1.1.3]{DAGX}. We recall some terminology from \cite{DAGX}:

\begin{definition}\label{smol}
A map $A \to B$ in $\mathscr A$ is called \emph{elementary} if it fits into a pullback square of the form \[\begin{tikzcd}
	A & {B} \\
	k & {\mrm{sqz}_k(k[m])}
	\arrow[from=1-1, to=1-2]
	\arrow[from=1-1, to=2-1]
	\arrow["\lrcorner"{anchor=center, pos=0.125}, draw=none, from=1-1, to=2-2]
	\arrow[from=1-2, to=2-2]
	\arrow[from=2-1, to=2-2]
\end{tikzcd}\] 
for some $m \ge 1$, where the bottom map is the `zero' map. More generally, a map $A \to B$ is \emph{small} if it can be written as a composite of elementary morphisms. An object $A$ is called small if the map $A \to k$ is small, and we denote the subcategory of small objects by $\mathscr A^{\mrm{art}}$.
\end{definition}

\begin{definition}\label{fmp-def}
We call a functor $F \colon \mathscr A^{\mrm{art}} \to \mathcal S$ a formal moduli problem if the following conditions hold:
\begin{enumerate}
	\item $F(k)$ is contractible.
	\item $F$ preserves pullbacks of cospans $A_0 \to A \leftarrow A_1$ of small maps.
\end{enumerate}
\end{definition}

We denote the corresponding full subcategories of $\Fun(\mathscr A^{\mrm{art}}, \mathcal S)$ by $\mrm{Moduli}^{\mathbf E_n}_{/k}$ and analogous variants.

\begin{remark} There are forgetful functors 
\[\mrm{Moduli}^{\mathbf E_0}_{\Lambda} \to \mrm{Moduli}^{\mathbf E_1}_{\Lambda} \to \dots \to \mrm{Moduli}^{\mathbf E_\infty}_{\Lambda} \to \mrm{Moduli}^{\mrm{an}}_{\Lambda};\]
at least heuristically, we can think of partition Lie algebras which arise from $\mrm{Moduli}^{\mathbf E_m}_\Lambda$ for small $m$ as being more trivial. For this reason, we shall at least initially consider categorical deformations over more general $\mathbf E_2$-algebras, even if the period map only concerns the underlying derived formal moduli problem.
\end{remark}

In many cases, we can make the descriptions much more explicit:

\begin{remark} An object of $\CAlg^{\an}_{\Lambda \sslash k}$, $\CAlg^{\an}_{/k}$, $\mrm{Alg}^{\mathbf E_n}_{\Lambda \sslash k}$ ($n \ge 1$), or $\mrm{Alg}^{\mathbf E_n}_{/k}$ ($n \ge 2)$ is small if and only if the underlying object $A$ in $\mrm{Alg}^{\mathbf E_1}_{/k}$ satisfies the following conditions:
  \begin{enumerate}
  \item $A$ is connective.
    \item The map $A \to k$ induces an equivalence $\pi_0(A)/J \xrightarrow{\simeq} k$, where $J$ denotes the Jacobson radical of $\pi_0(A)$.
  \item The (ordinary) ring $\pi_* A = \bigoplus_n \pi_n(A)$ is artinian.
  \end{enumerate}
This follows by adapting the proofs of \cite[Proposition 4.5.1]{DAGX} and \cite[Proposition 6.1.2]{DAGXII}, but we do not need this.

Similarly, small morphisms are surjective on $\pi_0$ (in fact for any $\mathscr A$ considered in this section), and the converse (for the above cases) holds by adapting \cite[Proposition 4.5.3]{DAGX} and \cite[Proposition 6.1.4]{DAGXII}; cf.\ also \cite[\S 3.1]{BT}. 

The restrictions on $\mathscr A$ above ensure that each $J^i/J^{i+1}$ is a symmetric $A$-module; this is not necessarily true otherwise.
\end{remark}

Finally, to fix notation, we recall that formal moduli problems have underlying classical deformation functors:
\begin{definition} Given a formal moduli problem $F \colon \mrm{CAlg}^{\mrm{an}, \mrm{art}}_{\Lambda \sslash k} \to \mathcal S$, the underlying classical deformation functor is the composite $F^{\mrm{cl}} \coloneq \pi_0 F \colon \mrm{CRing}^{\mrm{art}}_{\Lambda} \to \mrm{Set}$, where $\mrm{CRing}^{\mrm{art}}_{\Lambda} \subset \mrm{CAlg}^{\mrm{an},\mrm{art}}_{\Lambda}$ denotes the category of (discrete) Artinian local $\Lambda$-algebras with residue field $k$. 
\end{definition}

\subsection{Deforming schemes}

To fix notation, we recall the construction of the deformation functor for a derived scheme, closely following \cite[\S 3.6]{BT}.

\begin{construction}\label{scheme-construction}
Pulling back the source map $\Fun(\Delta^1, \mrm{DSch}^{\mrm{op}}) \longrightarrow \mrm{DSch}^{\mrm{op}}$ along $\mrm{Spec}$ produces a cocartesian fibration \[p\colon\mrm{CAlg}^{\mrm{an}}\mrm{DSch} \coloneq \mrm{CAlg}^{\mrm{an}} \times_{\mrm{DSch}^{\op}} \Fun(\Delta^1, \mrm{DSch}^{\op}) \longrightarrow \mrm{CAlg}^{\mrm{an}}.\]
The restriction of $p$ to the wide subcategory $\mrm{CAlg}^{\mrm{an}}\mrm{DSch}^{\mrm{cocart}}$ on the $p$-cocartesian edges is a left fibration. Given a derived $k$-scheme $X$, taking slices produces a left fibration \[\mrm{Defor}[X]\coloneq\mrm{CAlg}^{\mrm{an}}\mrm{DSch}^{\mrm{cocart}}_{/(k,X)} \to \mrm{CAlg}^{\mrm{an}}_{/k}, \] which corresponds to a functor $\mrm{Def}^+_{X}\colon \mrm{CAlg}^{\mrm{an}}_{/k} \longrightarrow \widehat{\mathcal S}$ via straightening (i.e.\ the $\infty$-categorical Grothendieck construction considered in \cite[Theorem 3.2.0.1, Proposition 3.3.2.5]{HTT}). We denote the restriction of $\mrm{Def}_X^{+}$ to $\mrm{CAlg}^{\mrm{art}}_{/k}$ (or $\CAlg^{\an}_{\Lambda}$) by $\Def_{X}$ (resp.\ $\mrm{Def}_{X, \Lambda}$).
\end{construction}

By \cite[Cor. 3.32]{BT}, $\mrm{Def}_{X}$ is valued in essentially small spaces, so it may be regarded as a functor into $\mathcal S$, and it is moreover a formal moduli problem.

\begin{remark}\label{cocart1} \begin{enumerate} \item An edge $e$ of $\mrm{CAlg}^{\mrm{an}}\mrm{DSch}$ is the data of a square \[\begin{tikzcd}
	Y' & {Y} \\
	{\Spec A'} & {\Spec A;}
	\arrow[from=1-1, to=1-2]
	\arrow[from=1-1, to=2-1]
	\arrow[from=1-2, to=2-2]
	\arrow[from=2-1, to=2-2]
      \end{tikzcd}\]
this square is cartesian if and only if $e$ is $p$-cocartesian. A point of $\mrm{Def}^+_X(A)$ is hence the data of a derived scheme $Y/A$ equipped with an equivalence $Y \times_A k \xrightarrow{\simeq} X$.
\item Since the fibre product above is formed in \emph{derived} schemes, $\mathrm{Def}_{X}(A)$ identifies with the usual $(1$-)groupoid of \emph{flat} lifts of $X$ whenever $A$ is discrete; see \cite[Corollary 3.35]{BT}.
\end{enumerate}
\end{remark}

\subsubsection{Tangent fibre}\label{looped-def}
It is not difficult to compute the tangent fibre of $\mrm{Def}_X$; we recall some of the details from \cite[\S 2-3]{BT}, which will be used later. Some notation is in order:
\begin{notation} Given a connective quasicoherent sheaf $\mathscr E \in \QC(X)_{\ge 0}$, we write $X^{\mathscr E}$ for the derived scheme $(|X|, \mathcal O_X \oplus \mathscr E)$. By abuse of notation, if $V \in \Mod_{k, \ge 0}$ and $X$ is a derived $k$-scheme, we sometimes write $X^V$ for $X^{\mathcal O_X \otimes_k V}$, where $\mathcal O_X \otimes_k V$ refers to the tensoring of $\QC(X)$ over $\Mod_k$.
\end{notation}
\begin{remark} If $A = \mrm{sqz}_k(V)$ for some $V \in \Mod_{k, \ge 0}$, then we may identify $X^V$ with the base change $X_A$, since $\mrm{sqz}_{\mathcal O_X}(\mathcal O_X \otimes_k V) \simeq \mathcal O_X \otimes_k \mrm{sqz}_k(V)$.
\end{remark}
For $V \in \Mod_{k, \ge 0}^{\omega}$ a connective perfect $k$-module, and $A = \mrm{sqz}_k(V)$, there are equivalences \begin{align*} 
\Omega \mrm{Def}_X(\mrm{sqz}_k(V)) &\simeq \mrm{Map}_{/A}(X_A, X_A) \times_{\mrm{Map}_{/k}(X, X)} \{\mrm{id}\} \\ 
&\simeq \mrm{Map}_{/k}(X_A, X) \times_{\mrm{Map}_{\mrm{DSch}_{k}}(X, X)} \{\mrm{id}\} \\
&\simeq \Map_{{X \sslash k}}(X_A, X) \\ 
&= \mrm{Map}_{X\sslash k}(X^V, X),
\end{align*}
where we have taken loops at the trivial deformation $X^V = X_A$.

For any $\mathscr E \in \QC(X)_{\ge 0}$, we may identify 
\[\mrm{Map}_{X\sslash k}(X^\mathscr E, X) \simeq \mrm{Map}_{k\sslash\mathcal O_X}(\mathcal O_X, \mathcal O_X \oplus \mathscr E) \simeq \mrm{Map}_{\QC(X)}(\mrm{L}_{X/k}, \mathscr E).\] Using the fact that $\mrm{Def}_X$ is a formal moduli problem, we then obtain equivalences \[ \mrm{Def}_X(\mrm{sqz}_k(V)) \simeq \Omega \mrm{Def}_X(\mrm{sqz}_k(V[1])) \simeq \mrm{Map}_{\QC(X)}(\mrm{L}_{X/k}, \mathcal O_X \otimes_k V[1]),\]
from which one readily computes the tangent fibre as $\mathrm{R}\Gamma(X,\mathrm{T}_{X})[1]$.

\subsubsection{Deformations of qcqs schemes}

It is convenient to assume that the categories of quasicoherent sheaves appearing are compactly generated; this happens in particular for $\QC(X)$ when $X$ is qcqs. This propagates to infinitesimal deformations:


\begin{proposition}\label{qcqs-propagation}
Let $A$ be an object of $\mrm{CAlg}^{\mrm{an, art}}_{/k}$, $Y$ be a derived $A$-scheme, and $X = Y \times_{A} k$. Then $Y$ is qcqs if and only if $X$ is. Hence if $X$ is qcqs, then $\mrm{Def}_X^{\mrm{qcqs}} \hookrightarrow \mrm{Def}_X$ is an equivalence.
\end{proposition}

\begin{proof}
Write $Z_{\mrm{cl}}$ for the underlying (ordinary) scheme of a derived scheme $Z$. Recall that in this situation, $Y_{\mrm{cl}} \to \Spec(\pi_0(A)) = \Spec(A)_{\mrm{cl}}$ is a flat lift of $X_{\mrm{cl}}$ (\cref{cocart1}), and, moreover, a map of derived schemes is qcqs if and only if the underlying map of (ordinary) schemes is. We are therefore reduced to proving the corresponding result about flat deformations of ordinary schemes, which follows from e.g.\ \cite[Lemma 09ZV (3) and (6)]{stacks-project}.
\end{proof}

\subsection{Deformations of linear $\infty$-categories}\label{catdefs}

\subsubsection{Recollections on linear/tensored $\infty$-categories}\label{lin-recoll}

We will use the formalism of linear (or \emph{tensored}) $\infty$-categories as described in \cite[Appendix D]{SAG}. Fix a compactly generated, presentably $\mathbf E_2$-monoidal stable $\infty$-category $\mathscr E$ and write $\mrm{Cat}^{\mrm{perf}}$ for the $\infty$-category of small idempotent complete stable $\infty$-categories and exact functors.

\begin{notation} We write
 \begin{enumerate}
   \item  $\PrL_{\mathscr E} \coloneq \mathrm{LMod}_{\mathscr E}(\mathrm{Pr}^{\mathrm{L}})$ for the $\infty$-category of left modules over $\mathscr E$ in $\PrL$, whose objects we call (presentably left) $\mathscr E$-linear $\infty$-categories, and $\Fun^{\mrm{L}}_{\mathscr E}(\mathscr M, \mathscr N)$ for the $\infty$-category of colimit-preserving $\mathscr E$-linear functors between two $\mathscr E$-linear $\infty$-categories $\mathscr M, \mathscr N$. 
   \item $\mrm{Cat}^{\mrm{perf}}_{\mathscr E} \coloneq \mrm{LMod}_{\mathscr E^\omega}(\mrm{Cat}^{\mrm{perf}})$ for the $\infty$-category of left modules over $\mathscr E^\omega \in \mrm{Alg}_{\mathbf E_2}(\mrm{Cat}^{\mrm{perf}})$.
   \end{enumerate}
 \end{notation}

Via $\mrm{Ind}$-completion, $\mrm{Cat}^{\mrm{perf}}_{\mathscr E^\omega}$ may be identified with the subcategory $\mrm{Pr}^{\mrm{L}, \omega}_{\mathscr E}$ of $\PrL_{\mathscr E}$ on the compactly generated $\mathscr E$-linear $\infty$-categories with compact-object-preserving functors; see \cite[\S 4]{Higher-Traces} for discussion, especially [\citelink{Higher-Traces}{Ibid}, Proposition 4.9.]. 

In practice, we will mostly concern ourselves with $\mathscr E = \mathrm{LMod}_A$ for an $\mathbf E_2$-ring $A$ and simply write $\mathrm{Cat}^{\mrm{perf}}_{A}$ for $\mrm{Cat}^{\mrm{perf}}_{\mrm{LMod}_A}$ and $\mathrm{Pr}^{\mrm{L}}_A$ for $\PrL_{\mrm{LMod}_A^\omega}$, calling them (small/presentable) $A$-linear $\infty$-categories. Moreover, given a presentable $A$-linear $\infty$-category $\mathscr C$ and a map $A \to B$ of $\mathbf E_2$-rings, we write $B \otimes_A \mathscr C$ for the left $B$-linear $\infty$-category $\mrm{LMod}_B \otimes_{\mrm{LMod}_A} \mathscr C$.

Given a colimit-preserving $\mathscr E$-linear functor $F$ between presentable $\mathscr E$-linear $\infty$-categories, the right adjoint $F^R$ to the (underlying functor of) $F$ admits a canonical \emph{lax} $\mathscr E$-linear structure (see \cite[Remark 7.3.2.9]{HA}). If $F^R$ is actually $\mathscr E$-linear, then the unit of the adjunction is canonically $\mathscr E$-linear; in other words, $F \dashv F^R$ in the $(\infty,2)$-category of presentable $\mathscr E$-linear $\infty$-categories. 

\begin{remark}\label{adjoints} If $\mathscr E^\omega$ is rigid (that is, every object is dualisable) and $F \colon \mathscr M \to \mathscr N$ is a map in $\PrL_\mathscr E$ where $\mathscr M$ is compactly generated, then 
\begin{enumerate}
\item The lax $\mathscr E$-linear right adjoint $G$ is $\mathscr E$-linear if it preserves filtered colimits, and 
\item $G$ preserves filtered colimits if and only if $F$ preserves compact objects; \end{enumerate}
this provides an `intrinsic' way to check whether $F$ has an $\mathscr E$-linear right adjoint. \end{remark}

\begin{proof} We could not quite find a reference for these, but they are certainly at least folklore and readily follow from the literature. For (1), it suffices to show that the maps $e \otimes G(x) \to G(e \otimes x)$ are equivalences for all $e \in \mathscr E$ and $x \in \mathscr N$, but the class of such $e$ contains all compact objects by \cite[Proposition 4.9]{Higher-Traces} (since $\mathscr E^\omega$ is rigid) and is closed under filtered colimits by our assumption (and since $\otimes$ preserves colimits in each variable).  For (2), note that $G$ preserving filtered colimits is equivalent to $\mrm{Map}(a, \colim_i G(b_i)) \to \mrm{Map}(a, G(\colim_i b_i))$ being an equivalence for all \emph{compact} $a$ and filtered diagrams $(b_i)$, but (using compactness of $a)$, this is equivalent to each $Fa$ being compact. \end{proof} 
Later we will also make use of the following basic lemma:
\begin{proposition}\label{expadj} \emph{(e.g. \cite[Lemma D.2.5.]{heyermann})}

Let $F\colon \mathscr M \to \mathscr N$ be an $\mathscr E$-linear functor with $\mathscr E$-linear right adjoint $G$ and unit $\eta$. Then for any other $\mathscr P  \in \PrL_{\mathscr E}$, the composition functor $F_*\colon \Fun_\mathscr E^{\mrm{L}}(\mathscr P, \mathscr M) \to \Fun_\mathscr E^{\mrm{L}}(\mathscr P,\mathscr N)$ is left adjoint to $G_*$ with unit $\eta_*$.  
\end{proposition}

\subsubsection{The deformation functor}

In this subsection, we recall the deformation functor governing deformations of a presentably $k$-linear $\infty$-category. This is a very mild variant of \cite[\S 5.3]{DAGX} adapted to mixed characteristic situations.

\begin{construction}\label{cat-construction} Pulling back the cocartesian fibration $\mrm{LMod}(\PrL) \longrightarrow \mrm{Alg}(\PrL)$ along the functor $\mrm{LMod}_{(-)} \colon \mrm{Alg}^{\mathbf E_2} \to \mrm{Alg}(\PrL)$, we obtain a cocartesian fibration \[q\colon\mrm{LinCat}^{\mathbf E_2}\coloneq \mathrm{Alg}^{\mathbf E_2} \times_{\mathrm{Alg}(\mathrm{Pr}^{\mrm{L}})}\mathrm{LMod}(\mathrm{Pr}^{\mrm{L}}) \longrightarrow \mrm{Alg}^{\mathbf E_2}.\] Objects in the source may be identified with pairs $(A \in \mrm{Alg}^{\mathbf E_2}, \mathscr C \in \PrL_A)$. Passing to the wide subcategory on $q$-cocartesian edges, we obtain a left fibration $\mrm{LinCat}^{\mathbf E_2, \mathrm{cocart}} \to\mathrm{Alg}^{\mathbf E_2}$.

Now let $k$ be a field and $\mathscr C$ a $k$-linear $\infty$-category. By \cite[Proposition 4.3.7.2]{kerodon}, taking slices produces a left fibration $\mrm{Defor}[\mathscr C]\coloneq \mrm{LinCat}^{\mathrm{cocart}}_{/(k, \mathscr C)} \longrightarrow \mathrm{Alg}^{\mathbf E_2}_{/k}$, classified by a functor $\mrm{CatDef}_{\mathscr C}^{+} \colon \mathrm{Alg}^{\mathbf E_2}_{/k} \to \widehat{\mathcal S}$. We write $\mrm{CatDef}_{\mathscr C, \Lambda}^{\mathbf E_2}$ for the further restriction to $\mrm{Alg}^{\mathbf{E}_2, \mrm{art}}_{\Lambda}$.
\end{construction}

\begin{variant}\label{an-cat}
As usual, one may restrict further to animated rings to produce functors $\mrm{CatDef}_{\mathscr C}^{+}$ and $\mrm{CatDef}_{\mathscr C}$; this corresponds to considering the pullback \[q^{\mrm{an}} \colon \mrm{LinCat}^{\mrm{an}} \coloneq \mrm{CAlg}^{\mrm{an}} \times_{\mrm{Alg}(\PrL)} \mrm{LMod}(\PrL) \longrightarrow \mrm{CAlg}^{\mrm{an}}\] of $q$ along $\mrm{CAlg}^{\mrm{an}} \to \mrm{Alg}^{\mathbf E_2}$. We write $\mrm{Defor}[\mathscr C]^{\mrm{an}}$ for the corresponding variant of $\mrm{Defor}[\mathscr C]$.
\end{variant}

\begin{remark} \label{cocart2}
\begin{enumerate}
\item $\mrm{CatDef}_{\mathscr C, k}^{\mathbf E_2}$ recovers Lurie's definition in \cite[Construction 16.6.1.2]{SAG}.
\item An edge $e\colon (A, \mathscr C) \to (B, \mathscr D)$ in $\mrm{LinCat}^{\mathbf E_2}$ or $\mrm{LinCat}^{\mrm{an}}$ is the datum of a pair $(f, \phi)$ of a map $f\colon A \to B$ and a $B$-linear map $\phi\colon B \otimes_{A} \mathscr C \to \mathscr D$. The edge $e$ is $q$-cocartesian if and only if $\phi$ is an equivalence (by \cite[Corollary 4.2.3.2]{HA}); in particular, we may identify $\mrm{CatDef}^+_{\mathscr C}(A)$ with the space of pairs $(\mathscr D, \psi \colon \mathscr D \otimes_A k \xrightarrow{\simeq} \mathscr C)$. 
\end{enumerate}
\end{remark}

\begin{variant} If $\mathscr C$ is compactly generated, one can also consider the subfunctors $\mrm{CatDef}_{\mathscr C}^{\mathbf E_2, \omega} \subset \mrm{CatDef}^{\mathbf E_2}_{\mathscr C}$ (and their variants) on those deformations which are compactly generated. 
\end{variant}

\subsection{Cohesiveness of $\mrm{CatDef}_{\mathscr C}$}

Unfortunately, $\mrm{CatDef}^{\mathbf E_2}_{\mathscr C}$ and its variants are not quite formal moduli problems in general, but this failure can be quantified. We recall from \cite[Definition 5.1.5]{DAGX} that a functor $F:\mathscr A^{\mrm{art}} \to \mathcal S$ (with $\mathscr A$ one of the usual deformation contexts of \cref{a-scr}) is called an \emph{$n$-proximate} formal moduli problem for some integer $n \ge 0$ if $F$ satisfies the conditions of \cref{fmp-def}, except that for a cospan of small maps $A \to B \leftarrow C$, the comparison map $F(A \times_B C) \to F(A) \times_{F(B)} F(C)$ is only required to be $(n-2)$-truncated.

We first state a general theorem due to Lurie. Here we recall that a stable $\infty$-category $\mathscr  C$ is called \emph{tamely compactly generated} if it is compactly generated and for each $c, d \in \mathscr C^\omega$, the mapping spectrum $\mrm{map}_{\mathscr C}(c,d)$ is almost connective.

\begin{theorem} \label{2prox} \emph{(Lurie)} Let $\mathscr C \in \PrL_k$.

\begin{enumerate} \item For every pullback diagram
  \[\begin{tikzcd} {B} \arrow[r] \arrow[d] & {B_0} \arrow[d] \\ {B_1} \arrow[r] & {B_{01}}
    \end{tikzcd}\]
  in $\mrm{Alg}^{\mathbf E_2}_{/k}$, the induced map $\mrm{CatDef}_{\mathscr C}^{+}(B) \to \mrm{CatDef}_{\mathscr C}^{+}(B_0) \times_{\mrm{CatDef}_{\mathscr C}^{+}(B_{01})} \mrm{CatDef}_{\mathscr C}^{+}(B_1)$ is $0$-truncated. In particular, $\mrm{CatDef}^{\mathbf E_2}_{\mathscr C} \colon \mrm{Alg}^{\mathbf E_2, \mrm{art}}_{/k} \to \widehat{\mathcal S}$ is a (large) $2$-proximate formal moduli problem.
\item If $\mathscr C$ is compactly generated, then $\mrm{CatDef}^{\mathbf E_2, \omega} \colon \mrm{Alg}^{\mathbf{E}_2, \mrm{art}}_{/k} \to \widehat{\mathcal S}$ is valued in essentially small spaces and, considered as a functor into $\mathcal S$, is a $2$-proximate formal moduli problem.
\item \label{tcg} If $\mathscr C$ is \emph{tamely} compactly generated, then $\mrm{CatDef}^{\mathbf E_2, \omega} \colon \mrm{Alg}^{\mathbf E_2, \mrm{art}}_{/k} \to \mathcal S$ is a $1$-proximate formal moduli problem.
\end{enumerate}
\end{theorem}

\begin{proof} Let us explain why this follows from results of \cite{SAG}. Variants of (1) and (2) with $\mrm{Alg}^{\mathbf E_2}_{/k}$ replaced by $\mrm{Alg}^{\mathbf E_2}_{k \sslash k}$ are proven as Propositions 16.6.2.1 and 16.6.6.2 respectively, but the proofs never use the $k$-algebra structure. Part (3) is a slight strengthening of Theorem 16.6.7.3, which Lurie deduces from Propositions 16.6.9.1 and 16.6.9.2. The proof of the former manifestly never uses the $k$-algebra structure, whilst the latter goes through by noting that the dévissage reducing to square-zero extensions still goes through by definition of $\mrm{Alg}^{\mathbf E_2, \mrm{art}}_{/k}$. 
\end{proof}

Finally, let us comment on approximation by formal moduli problems. In general, for any deformation context $\mathscr A$, the inclusion $\mrm{Moduli}_{\mathscr A} \hookrightarrow \Fun(\mathscr A, \mathcal S)$ admits a left adjoint $(-)^{\land}$. If $\mathscr A$ admits a deformation theory, then for any $n$-proximate formal moduli problem $F \colon \mathscr A \to \mathcal S$, the unit map $F \to F^{\land}$ is $(n-2)$-truncated \cite[Theorem 5.1.9]{DAG}. A proof that all the deformation contexts considered in \cref{a-scr} admit a deformation theory does not appear in the literature, but by \cite[Theorem 6.25 (2)]{BT}, $\mrm{CAlg}^{\mrm{an}}_{\Lambda \sslash k}$ does. We therefore have the following:
\begin{corollary} \label{completion} If $\mathscr C$ is tamely compactly generated, then the map $\mrm{CatDef}_{\mathscr C, \Lambda}^{\omega} \to \mrm{CatDef}_{\mathscr C, \Lambda}^{\land}$  is $(-1)$-truncated. 
\end{corollary}
 
\begin{remark}
Recent work of Iwanari shows (at least when $\Lambda = k$) that one obtains a formal moduli problem if one considers deformations equipped with left-complete $t$-structures \cite{iwanari-t-structures}. Our period map could be adapted to factor through this formal moduli problem; cf.\ [\citelink{iwanari-t-structures}{Ibid}, Example 4.4]. However, this does not make a difference on tangent fibres, and neglecting $t$-structures allows for stronger conclusions regarding derived equivalences between schemes.
\end{remark}

\subsection{$\QC(X)$ and Hochschild (co)homology}\label{hochschild-recoll} Let us describe how the discussion of the previous section applies to $\infty$-categories of quasicoherent sheaves on derived schemes.

We first discuss the linear structure. Taking quasicoherent sheaves refines to a functor \[\QC \colon \mrm{DSch}^{\mrm{op}} \to \mrm{CAlg}(\PrL).\] By \cite[\S 3.4.1.3]{HA}, there is an equivalence \[\mrm{CAlg}(\PrL)_{\mathscr C/} \simeq \mrm{CAlg}(\Mod_{\mathscr C}(\PrL))\] for each $\mathscr C \in \mrm{CAlg}(\PrL)$, which informally sends a map $\mathscr C \to \mathscr D$ of symmetric monoidal $\infty$-categories to $\mathscr D$ viewed as a $\mathscr C$-linear symmetric monoidal $\infty$-category. Therefore, for any derived scheme $Y$, we obtain a lift \[\QC \colon \big(\mrm{DSch}_{/Y}\big)^{\op} \to \Mod_{\QC(Y)}(\PrL).\] In particular, $\QC(Z)$ has a canonical $A$-linear structure for any ring $A$ and derived $A$-scheme $Z$.

\subsubsection{Quasicoherent sheaves on qcqs derived schemes }

It is convenient at some points to ensure that the $\infty$-category of quasicoherent sheaves on a derived scheme is compactly generated. We highlight a sufficient condition for this and some consequences; this is all well-known, but we could not find it explicitly stated in the literature in this generality.

\begin{proposition}\label{qcqs-cats} 
Let $Y,Z,W, S$ be derived schemes, and $W \to S \leftarrow Z$ be any cospan. 
\begin{enumerate}
	\item \label{qcqs-tcg} If $Y$ is qcqs, then $\QC(Y)$ is tamely compactly generated, and an object $\mathscr F$ is compact if and only if it is perfect. 
	\item \label{tensor} If $S$ is qcqs and $\QC(W)$ is dualisable (e.g.\ if it is compactly generated), then the canonical comparison map 
	\[ \QC(W) \otimes_{\QC(S)} \QC(Z) \longrightarrow \QC(W \times_{S} Z)\] is an equivalence. 
	\item \label{fm}  Let $\pi_1, \pi_2$ be the two projections out of $W \times_{S} Z$. Under the assumptions of \cref{tensor}, the `Fourier--Mukai' assignment \[\mathscr F \mapsto (\pi_{2})_*(\pi_1^*(-)  \otimes \mathscr F)\] induces an equivalence 
	\[ \QC({W} \times_{S}Z) \xlongrightarrow{\simeq} \mrm{Fun}^{\mrm{L}}_{\QC(S)}(\QC(W), \QC(Z)).\] 
\end{enumerate}
\end{proposition}

We will show how to deduce this from various results in the literature.
\begin{proof}
\hyperref[qcqs-tcg]{Part \ref*{qcqs-tcg}} follows from \cite[Proposition 9.6.1.1, Corollary 9.6.3.2]{SAG} if we consider $Y$ as a spectral algebraic space; the compact generation is tame since perfect complexes on a quasicompact derived scheme are bounded. \hyperref[tensor]{Part \ref*{tensor}} follows from \cite[Proposition 3.5.3]{GaitsgoryRozenblyum}, which applies in virtue of \hyperref[qcqs-tcg]{Part \ref*{qcqs-tcg}}. Finally, \hyperref[fm]{(\ref*{fm})}  follows by self-duality of $\QC(W)$ as a $\QC(S)$-module (e.g. \cite[Corollary 9.4.2.2]{SAG}), and because \[\mrm{Fun}^{\mrm{L}}_{\QC(S)}(\QC(W),\QC(Z)) \simeq \QC(W)^\vee \otimes_{\QC(S)} \QC(Z),\] the dual being taken in $\QC(S)$-modules; cf.\ the proofs of \cite[Corollary 4.10.]{BZFN} and \cite[Theorem 8.9.]{toenfm}.
\end{proof}

\begin{remark} \label{tame} Combined with \cref{qcqs-cats} (\ref{qcqs-tcg}), \hyperref[tcg]{Theorem \ref*{2prox} (\ref*{tcg})} implies that $\mrm{CatDef}_{\QC(X)}^\omega$ is a $1$-proximate formal moduli problem for any qcqs derived $k$-scheme $X$.
\end{remark}
\subsubsection{The looped moduli problem and Hochschild cohomology}
It follows from \cite[Theorem 5.3.16]{DAGX} that the tangent fibre of the $\mathbf E_2$-formal moduli problem $\mrm{CatDef}_{\mathscr C, k \sslash k}^{\land}$ is given by the shifted Hochschild cochains $\mrm{HH}^\bullet(\mathscr C/k)[2]$, where $\mrm{HH}^\bullet(\mathscr C/k)$ denotes the mapping $k$-module
\[\mrm{HH}^\bullet(\mathscr C/k) = \mrm{map}_{\mrm{Fun}^{\mrm{L}}_k(\mathscr C, \mathscr C)}(\mrm{id}, \mrm{id}).\]
A key part of Lurie's proof computes $\Omega^2 \mrm{CatDef}_{\mathscr C}$ at square-zero extensions. In the situation where $\mathscr C = \QC(X)$ for a qcqs derived $k$-scheme $X$, it is possible to deloop this to a computation of $\Omega \mrm{CatDef}^\omega_{\QC(X)}(\mrm{sqz}_k(V))$ for $V \in \Mod_{k, \ge 0}$; this will later ease comparison with $\Omega \mrm{Def}_X$.

\begin{notation} Given a map $y \to x$ in an $\infty$-category $\mathscr C$, we write \[\mrm{Map}'_{\mathscr C}(x,y) \coloneq \mrm{Map}_{\mathscr C}(x,y) \times_{\mrm{Map}_{\mathscr C}(x,x)} \{\mrm{id}\}\] for the space of sections of $y \to x$. Usually the map $y \to x$ will be clear from context, so there is no risk of confusion. 
\end{notation}

The first step towards computing $\Omega \mrm{CatDef}$ is the following:
\begin{proposition}\label{loops-catdef} Let $X/k$ be a qcqs derived scheme. There is an equivalence
\[\Omega \mrm{CatDef}_{\QC(X), \Lambda \sslash k}(\mathrm{sqz}_k(V)) \simeq \mrm{Map}'_k(\QC(X), \QC(X^V))\]
natural in $V \in \Mod_{k, \ge 0}$, where the implicit map $\QC(X^V) \to \QC(X)$ is pullback along $X \to X^V$.
\end{proposition}
\begin{proof}
Unwinding definitions and putting $A = \mathrm{sqz}_k(V)$, we have equivalences
\[ \Omega \mrm{CatDef}_{\QC(X)}(A) \simeq \mrm{Map}_A(\QC(X) \otimes_k A, \QC(X) \otimes_k A) \times_{\mrm{Map}_k(\QC(X), \QC(X))} \{\mrm{id}\},\] where $(-) \otimes_k A$ is the extension of scalars functor $\PrL_k \to \PrL_A$;  cf.\ the description in \cite[16.6.2]{SAG}. By adjunction,
\begin{align*}
\Omega \mrm{CatDef}_{\QC(X)}(A) &\simeq \mrm{Map}_k(\QC(X), \QC(X) \otimes_k A) \times_{\mrm{Map}_k(\QC(X), \QC(X))} \{\mrm{id}\}  \\ &= \mrm{Map}_k'(\QC(X), \QC(X) \otimes_k A) \\
&\simeq \mrm{Map}_k'(\QC(X), \QC(X \times_k A)) \\ &= \mrm{Map}_k'(\QC(X), \QC(X^V)),
\end{align*} where we have used \cref{qcqs-cats} (\ref*{tensor}) to identify $\QC(X) \otimes_k A \simeq \QC(X \times_k A)$.
\end{proof}

\begin{remark}\label{fm2} Under the `Fourier--Mukai' equivalence $\Fun_k^{\mrm{L}}(\QC(X), \QC(Y)) \simeq \QC(X \times_k Y)$ of \hyperref[fm]{Theorem \ref*{qcqs-cats} (\ref*{fm})}, the diagonal $\Delta_* \mathcal O_X$ corresponds to the identity functor. Taking mapping $k$-module spectra leads to equivalences
\[ \mrm{HH}^\bullet(X/k) \simeq \mrm{map}_{\QC(X \times_k X)}(\Delta_* \mathcal O_X, \Delta_* \mathcal O_X) \simeq \mrm{map}_{\QC(X)}(\Delta^* \Delta_* \mathcal O_X, \mathcal O_X)\] which relate Hochschild cohomology to the Hochschild homology sheaf $\underline{\mrm{HH}}_\bullet \simeq \Delta^* \Delta_* \mathcal O_X$ (see \cref{relating}).
\end{remark}
\section{Atiyah classes and transformations}\label{hochschild-atiyah}

To better understand $\mrm{CatDef}_{\QC(X)}$, and eventually the commutative-to-noncommutative period map, we now wish to compute the spaces 
\[ \Omega \mrm{CatDef}_{\QC(X)}(\mrm{sqz}_k(V)) \simeq \mrm{Map}'_{k}(\QC(X), \QC(X^V))\] appearing in \cref{loops-catdef} in more geometric terms. For later applications, we will more generally study square-zero extensions $X^{\mathscr E}$ by \emph{any} connective quasicoherent sheaf $\mathscr E$ (rather than those of the form $\mathcal O_X \otimes_k V$). The aim of this section is to prove the following:

\begin{restatable}{theorem}{hkrat}
\label{HKR-atiyah} Let $X/k$ be a qcqs derived $k$-scheme. 
\begin{enumerate}
	\item \label{one} There are canonical equivalences of mapping spaces \[\mrm{Map}_k'(\QC(X), \QC(X^{\mathscr E})) \simeq \mrm{Map}_{\Fun_k^{\mrm{L}}(\QC(X), \QC(X))}(\mrm{id}, - \otimes \mathscr E[1]) \simeq \mrm{Map}_{\QC(X)}(\underline{\HH}_{X/k, \bullet}, \mathscr E[1]),\] natural in $\mathscr E \in \QC(X)_{\ge 0}$. 
	\item \label{two} When $\mathscr E = \mrm{L}_{X/k}$ is the cotangent complex, the composite equivalence of \hyperref[one]{\emph{(\ref*{one})}} sends the pullback $s_{\delta}^*$ along the map $s_\delta \colon X^{\mrm{L}_{X/k}} \to X$ (induced by the universal derivation $\mathcal O_X \to \mrm{L}_{X/k}$) to the \emph{HKR} map $\underline{\mrm{HH}}_{\bullet} \longrightarrow \mrm{L}_{X/k}[1]$.
\end{enumerate}
\end{restatable}

The proof will make use of generalised Atiyah classes, which we connect to Hochschild homology by generalising results of Căldăraru and Buchweitz--Flenner (\cref{HKR-Atiyah-Lemma}).

\begin{notation}
Given a square-zero extension $X^\mathscr E$ of a derived scheme $X$ by some $\mathscr E \in \QC(X)_{\ge 0}$, we write $i$ for the inclusion map $X \to X^\mathscr E$. Given a right inverse $f \colon \mathcal O_X \to \mrm{sqz}_{\mathcal O_X}(\mathscr E)$ to the projection, we write $s_f$ for the induced map $X^\mathscr E \to X$.
\end{notation}

\subsection{The derived Atiyah transformation}
In this subsection, we describe what we shall call the \emph{derived} Atiyah transformation.

First, let us briefly recall Lurie's construction of Atiyah classes in \cite[\S 19.2.2]{SAG}. Fixing a spectral Deligne--Mumford stack $\mathsf Y$ and some $\mathscr F \in \QC(\mathsf Y)_{\gg -\infty}$, Lurie constructs an equivalence \[\QC(\mathsf Y^{\mathscr E}) \times_{\QC(\mathsf{Y})} \{\mathscr F\} \simeq \mrm{Map}_{\QC(\mathsf Y)}(\mathscr F, \mathscr F \otimes \mathscr{E}[1])\] functorial in $\mathscr E \in \QC(\mathsf Y)_{\ge 0}$ \cite[Proposition 19.2.2.2]{SAG}; this in particular classifies lifts of $\mathscr F$ to ${\mathsf Y}^{\mathscr E}$ modulo equivalence in terms of $\mrm{Ext}^1$-classes. Denoting the $\mathbf E_\infty$ (or \emph{topological}) cotangent complex of $\mathsf Y$ by $\mrm{L}_{\mathsf Y}^{\mathbf E_\infty}$, the universal derivation $\mathcal O_{\mathsf Y} \to \mrm{L}_{\mathsf Y}^{\mathbf E_\infty}$ induces a right inverse $s_{\delta} \colon {\mathsf Y}^{\mrm{L}_{\mathsf Y}^{\mathbf{E}_\infty}} \to \mathsf Y$ of $\mathsf Y \to \mathsf Y^{\mrm{L}_{\mathsf Y}^{\mathbf E_\infty}}$. The pullback $s_{\delta}^* \mathscr F$ therefore provides a lift of $\mathscr F$ to ${\mathsf Y}^{\mrm{L}_{\mathsf Y}^{\mathbf E_\infty}}$, and hence corresponds to a class in $\mrm{Ext}^1(\mathscr F, \mathscr F \otimes \mrm{L}^{\mathbf E_\infty}_{\mathsf Y})$ which Lurie calls the \emph{Atiyah class} $\mrm{At}(\mathscr F)$ of $\mathscr F$ \cite[Construction 19.2.2.5]{SAG}. Let us call this Lurie's \emph{spectral} Atiyah class. This construction readily adapts to the setting of derived algebraic geometry to produce what we shall call Lurie's (derived) Atiyah classes.

In this subsection, we consider a variant of Lurie's construction which concerns `deforming all quasicoherent sheaves at once' (\cref{classify-sections}); this leads to a \emph{derived Atiyah transformation} $\mrm{id} \longrightarrow - \otimes \mrm{L}_{X/k}[1]$ which specialises pointwise to the derived Atiyah class. We first require two technical lemmata:
\begin{lemma}\label{cartesian} \emph{(Lurie)} Given a pushout of derived schemes \[\begin{tikzcd}
	{Y_{01}} & {Y_0} \\
	{Y_1} & Y
	\arrow["{{{q_0}}}"{inner sep=.8ex}, from=1-1, to=1-2]
	\arrow["{{{q_1}}}"'{inner sep=.8ex}, from=1-1, to=2-1]
	\arrow["{{{f_0}}}", from=1-2, to=2-2]
	\arrow["{{{f_1}}}"', from=2-1, to=2-2]
	\arrow["\lrcorner"{anchor=center, pos=0.125, rotate=180}, draw=none, from=2-2, to=1-1]
\end{tikzcd}\]
along closed immersions, the comparison  map \[\QC(Y) \longrightarrow \QC(Y_0) \times_{\QC(Y_{01})} \QC(Y_1)\] is fully faithful, and restricts to an equivalence $\Perf(Y) \xlongrightarrow{\simeq} \Perf(Y_0) \times_{\Perf(Y_{01})} \Perf(Y_1)$. 
\end{lemma}

\begin{proof} We deduce this from \cite[\S 16.2]{SAG}, and in particular Theorem 16.2.0.1, following Lurie's arguments. By Zariski descent for quasicoherent sheaves, the assertions are local on $Y$, so we may assume $Y = \Spec A$ for some animated ring $A$. In this situation, the four derived schemes involved are all affine, and the pushout of derived schemes corresponds to a pullback
\begin{equation}\label{pullback}
\begin{tikzcd}
	A & {A_0} \\
	{A_1} & {A_{01}}
	\arrow[from=1-1, to=1-2]
	\arrow[from=1-1, to=2-1]
	\arrow["\lrcorner"{anchor=center, pos=0.125}, draw=none, from=1-1, to=2-2]
	\arrow[from=1-2, to=2-2]
	\arrow[from=2-1, to=2-2]
\end{tikzcd}
\end{equation}
of animated rings where the maps are surjective on $\pi_0$. By \cite[Theorem 16.2.0.2]{SAG}, the induced comparison map \[\Mod_{A} \longrightarrow \Mod_{A_0} \times_{\Mod_{A_{01}}} \Mod_{A_1}\] is fully faithful, proving the first assertion, and restricts to an equivalence \[\Mod_{A, \ge n} \xlongrightarrow{\simeq} \Mod_{A_0, \ge n} \times_{\Mod_{A_{01}, \ge n}} \Mod_{A_1, \ge n}\] for any $n \in \mathbf Z$. By [\citelink{SAG}{Ibid}, Proposition 16.2.3.1], an object of $\Mod_A$ is perfect if and only if its images in both $\Mod_{A_0}$ and $\Mod_{A_1}$ are. Since perfect modules over connective $\mathbf E_1$-rings are bounded below, we see that the comparison map restricts to an equivalence 
\[\Mod_{A}^\omega \xlongrightarrow{\simeq} \Mod_{A_0}^\omega \times_{\Mod_{A_{01}}^\omega} \Mod_{A_1}^{\omega},\] which implies the second claim of the Lemma.
\end{proof}

\begin{lemma}\label{compacts} Let $X$ be a qcqs derived $k$-scheme. Then for any $\mathscr E \in \QC(X)_{\ge 0}$,
\begin{enumerate}
	\item The functor $i^* \colon \QC(X^\mathscr E) \to \QC(X)$ reflects compact objects.
	\item The inclusion $\mrm{Map}_{k,\omega}'(\QC(X),\QC(X^{\mathscr E})) \subset \mrm{Map}_{k}'(\QC(X), \QC(X^{\mathscr E}))$ is an equivalence,  where these are spaces of sections as in \cref{loops-catdef}.
\end{enumerate}
\end{lemma}
\begin{proof}
For (1), write $X^{\mathscr E} = X \coprod_{X^{\mathscr E[1]}} X$ and use the result (used in the proof of \cref{cartesian}) that a quasicoherent sheaf on $X^{\mathscr E}$ is perfect if and only if its restrictions to $X$ and $X^{\mathscr E[1]}$ are.  (1) now implies that any section of $i^*$ preserves compact objects, so (2) follows immediately.
\end{proof}

We can now construct the first equivalence of \cref{HKR-atiyah} (1); cf.\ \cite[Proposition 19.2.2.2]{SAG}.
\begin{proposition}\label{classify-sections} 
	
Let $X$ be a qcqs derived $k$-scheme. There is a canonical equivalence
\[\mrm{Map}_k'(\QC(X),\QC(X^{\mathscr E})) \xlongrightarrow{\simeq} \mrm{Map}_{\Fun_k^{\mrm{L}}(\QC(X), \QC(X))}(\mrm{id}, (-) \otimes \mathscr E[1]),\] natural in $\mathscr E \in \QC(X)_{\ge 0}$. 
\end{proposition}
\begin{proof}
There is a pushout diagram
\[\begin{tikzcd}
	{X^{\mathscr E[1]}} & X \\
	X & {X^{\mathscr E}}
	\arrow["u", from=1-1, to=1-2]
	\arrow["u"', from=1-1, to=2-1]
	\arrow["i", from=1-2, to=2-2]
	\arrow["i"', from=2-1, to=2-2]
	\arrow["\lrcorner"{anchor=center, pos=0.125, rotate=180}, draw=none, from=2-2, to=1-1]
\end{tikzcd}\] in $\mrm{DSch}_{k}$, induced by a corresponding pullback in $\mrm{Shv}(X, \mrm{CAlg}^{\mrm{an}}_k)$, whose formation is functorial in $\mathscr E$. Applying the functor $\QC(-)$, we obtain a square in $\mathrm{Pr}^{\mathrm{L},{\omega}}_k$ which is cartesian by \cref{cartesian}:
\[\begin{tikzcd}
	{\QC(X^{\mathscr E})} & {\QC(X)} \\
	{\QC(X)} & {\QC(X^{\mathscr E[1]})}
	\arrow["{i^*}", from=1-1, to=1-2]
	\arrow["{i^*}", from=1-1, to=2-1]
	\arrow["\lrcorner"{anchor=center, pos=0.125}, draw=none, from=1-1, to=2-2]
	\arrow["{u^*}", from=1-2, to=2-2]
	\arrow["{u^*}", from=2-1, to=2-2]
\end{tikzcd}\]
We now apply $\mrm{Map}_{k}(\QC(X),-)$ and take fibres over $\mrm{id}_{\QC(X)}$ to obtain a pullback of spaces
\[\begin{tikzcd}
	{\mrm{Map}'_k(\QC(X),\QC(X^{\mathscr E}))} & {\mrm{Map}'_k(\QC(X),\QC(X))} \\
	{\mrm{Map}_k'(\QC(X),\QC(X))} & {\mrm{Map}_k'(\QC(X), \QC(X^{\mathscr E[1]}))} 
	\arrow[from=1-1, to=1-2]
	\arrow[from=1-1, to=2-1]
	\arrow[from=1-2, to=2-2]
	\arrow[from=2-1, to=2-2]
\end{tikzcd}\]
Here we have used \cref{compacts} to replace $\mrm{Map}_{k, \omega}$ with $\mrm{Map}_k$.

As the space $\mrm{Map}'_k(\QC(X), \QC(X))$ is contractible, we obtain natural equivalences 
\begin{align*}
\mrm{Map}'_k(\QC(X), \QC(X^{\mathscr E})) &\simeq \Omega \big(\mrm{Map}_k'(\QC(X), \QC(X^{\mathscr E[1]})), \chi \big) \\
& \simeq \Omega \big(\mrm{Map}_k(\QC(X), \QC(X^{\mathscr E[1]})), u^*\big) \times_{\Omega(\mrm{Map}_k(\QC(X), \QC(X)), \mrm{id})} \{\mrm{id}\} \\
& \simeq \mrm{Map}_{\Fun_k^{\mrm{L}}(\QC(X), \QC(X^{\mathscr E[1]}))}(u^*, u^*) \times_{\mrm{Map}_{\Fun_k^{\mrm{L}}(\QC(X), \QC(X))}(\mrm{id}, \mrm{id})} \{ \mrm{id}\}
\end{align*} 
Here $\Omega(S,s)$ denotes the loop space of a space $S$ at $s$, and we write $\chi$ for the point $(u^*, v^* u^* \simeq \mrm{id})$ of $\mrm{Map}_k'(\QC(X),\QC(X^{\mathscr E[1]}))$ ($v$ being the embedding $X~ \to~X^{\mathscr E[1]}$). We have used the fact that $v^*$ is conservative to replace $\mrm{Map}_k$ by $\Fun_k^{\mrm{L}}$.

Applying \cref{expadj} to the adjunction between $u^*$ and $u_*$, this space is in turn equivalent to
\begin{align*} &\mrm{Map}_{\Fun_k^{\mrm{L}}(\QC(X), \QC(X))}(\mrm{id}, u_*u^*) \times_{\mrm{Map}_{\Fun_k^{\mrm{L}}(\QC(X), \QC(X))}(\mrm{id}, \mrm{id})} \{ \mrm{id}\} \\
\simeq & \mrm{Map}_{\Fun_k^{\mrm{L}}(\QC(X), \QC(X))}(\mrm{id}, \mrm{id} \oplus (-) \otimes \mathscr E[1]) \times_{\mrm{Map}_{\Fun_k^{\mrm{L}}(\QC(X), \QC(X))}(\mrm{id}, \mrm{id})} \{\mrm{id}\} \\
\simeq & \mrm{Map}_{\Fun_k^{\mrm{L}}(\QC(X), \QC(X))}(\mrm{id}, - \otimes \mathscr E[1]). \qedhere
\end{align*}
\end{proof}

\begin{remark} Unravelling definitions, the equivalence of \cref{classify-sections} admits the following description. Any object $\sigma \in \mrm{Map}_k(\QC(X), \QC(X^\mathscr E))$ determines an object \[(i^* \sigma, i^* \sigma, \psi) \in \mrm{Map}_k(\QC(X), \QC(X)) \times_{\mrm{Map}_k(\QC(X),\QC(X^{\mathscr E[1]}))} \mrm{Map}_k(\QC(X),\QC(X)),\] where $\psi: u^* i^* \sigma \simeq u^* i^* \sigma$ is an equivalence. An equivalence $i^* \sigma \simeq \mrm{id}$ witnessing $\sigma$ as a section of $i^*$ induces an equivalence \[(i^* \sigma, i^* \sigma, \psi) \simeq (\mrm{id}, \mrm{id}, \beth(\sigma))\]  for some $\beth(\sigma): u^* \simeq u^*$ (this notation $\beth(\sigma)$ is slightly abusive, as it depends on the choice of identification $i^* \sigma \simeq \mrm{id}$, but we will almost exclusively deal with sections equipped with a `witness').  The induced map $\mrm{id} \to - \otimes \mathscr E[1]$ is then the composite 
\[\mrm{id} \xrightarrow{\overline{\beth(\sigma)}} u_* u^* = \mrm{id} \oplus (-) \otimes \mathscr E[1] \longrightarrow - \otimes  \mathscr E[1],\] where $\overline{\beth(\sigma)}$ is adjoint to $\beth(\sigma)\colon u^* \simeq u^*$.
\end{remark}

\begin{construction}\label{derivedatiyah} Given a point of $\mrm{Map}(X^{\mathscr E}, X) \times_{\mrm{Map}(X, X)} \{\mrm{id}\}$, i.e.\ a map $r \colon X^{\mathscr E} \to X$ along with an equivalence $ri \simeq \mrm{id}$, taking the pullback provides a point $\mrm{Map}'_k(\QC(X), \QC(X^\mathscr E))$, and hence a map $\mrm{id} \to (-) \otimes \mathscr E[1]$ of $k$-linear functors from \cref{classify-sections}. In the universal case that $\mathscr E = \mrm{L}_{X/k}$ and we take $r = s_\delta$ induced by the universal derivation $\mrm{d} \colon \mathcal O_X \to \mrm{L}_{X/k}$, we obtain a natural transformation \[\mathrm{At}^{\mrm{der}} \colon \mrm{id} \longrightarrow (-) \otimes \mrm{L}_{X/k}[1],\] which we call the \emph{derived Atiyah transformation}.
\end{construction}

\begin{remark} 
By construction, the derived Atiyah transformation specialises pointwise to the derived Atiyah class: for each $\mathscr E \in \QC(X)_{\ge 0}$ and $\mathscr F \in \QC(X)$, there are commutative diagrams 
\[\begin{tikzcd}
	{\mrm{Map}_k(\QC(X), \QC(X^{\mathscr E})) \times_{\mrm{Map}_k(\QC(X), \QC(X))} \{\mrm{id}\}} & {\mrm{Map}_{\Fun_k^{\mrm{L}}(\QC(X), \QC(X))}(\mrm{id}, (-) \otimes \mathscr E[1])} \\
	{\QC(X^{\mathscr E}) \times_{\QC(X)}\{\mathscr F\} } & {\mrm{Map}(\mathscr F, \mathscr F \otimes \mathscr E[1])}
	\arrow["\simeq", from=1-1, to=1-2]
	\arrow[from=1-1, to=2-1]
	\arrow[from=1-2, to=2-2]
	\arrow[from=2-1, to=2-2]
\end{tikzcd}\] where the vertical maps are given by evaluation at $\mathscr F$. 
\end{remark}
 
\subsection{Another model for Atiyah transformations}\label{classical}

We shall now compare the above derived Atiyah transformation to more classical definitions.

\begin{construction}\label{classicalatiyah} Let $X$ be a qcqs derived $k$-scheme. There is a canonical fibre sequence \[i_* \mrm{L}_{X/k} \to \mrm{sqz}_{\mathcal O_X}(\mrm{L}_{X/k}) \to i_* \mathcal O_X\] in $\QC(X^{\mrm{L}_{X/k}})$, known as the sequence of principal parts. We can consider its cofibre $i_* \mathcal O_X \to i_* \mrm{L}_{X/k}[1]$. We can push this forward along the map $\Delta' = (s_0, s_\delta) \colon X^{\mrm{L}_{X/k}} \to X \times_k X$ and use the equivalence $\Delta \simeq \Delta' i$ to obtain a map $\Delta_* \mathcal O_X \to \Delta_* \mrm{L}_{X/k}[1]$. Interpreting this as a map of Fourier--Mukai kernels, this induces a canonical $k$-linear natural transformation $\mrm{id} \to - \otimes \mrm{L}_{X/k}[1]$. This is called the (classical) \emph{Atiyah transformation}.
\end{construction}
  
\begin{remark}
We can reformulate the definition of the Atiyah transformation as follows. The map $\mrm{sqz}_{\mathcal O_X}(\mrm{L}_{X/k}) \to \mathcal O_X$ appearing is simply the adjunction unit $\mrm{id} \to i_* i^*$ evaluated at $\mrm{sqz}_{\mathcal O_X}(\mrm{L}_{X/k})$, and the natural transformation above is equivalently the cofibre $\mrm{id} \to - \otimes \mrm{L}_{X/k}[1]$ of the canonical map $(s_0)_* s_{\delta}^* \xrightarrow{\mrm{unit}} {(s_{0})}_{*} i_* i^* s_{\delta}^* \simeq \mrm{id}$. 
\end{remark}

\begin{remark}(Some comparisons)

\begin{itemize} \item For each $\mathscr F \in \QC(X)$, the induced class in $\mrm{Ext}^1(\mathscr F, \mathscr F \otimes \mrm{L}_{X/k})$ is called the \emph{Atiyah class} $\mrm{At}(\mathscr F)$ of $\mathscr F$. By construction, this agrees with Illusie's construction in \cite[IV \S 2.3.]{Illusie} (when the latter is defined).
  \item Here $X^{\mrm{L}_{X/k}}$ is a replacement for the infinitesimal neighbourhood $X^{(2)}$ of the diagonal. If $X/k$ is smooth, then there is a canonical equivalence $X^{\mrm{L}_{X/k}} \simeq X^{(2)}$ under which $\Delta'$ corresponds to the usual inclusion $X^{(2)} \hookrightarrow X \times_k X$; see \cite[Lemma 7.8]{rienks}.
  \end{itemize}
\end{remark}

We shall verify the following:
\begin{proposition} \label{comparison}
The derived Atiyah transformation of \cref{derivedatiyah} is equivalent to the classical Atiyah transformation of \cref{classicalatiyah}.
\end{proposition} 

\begin{remark}\label{mundinger-remark} In private communication, J. Mundinger sketched an independent proof for the case of Atiyah classes. In characteristic $0$, essentially the same result appears as \cite[Proposition 3.40]{Hennion-Lie}. \end{remark}

Recall that, as in \cref{classify-sections}, given a point $(\sigma, \alpha \colon i^* \sigma \simeq \mrm{id}) \in \mrm{Map}'(\QC(X), \QC(X^{\mathscr E}))$, the image of $\sigma$ in \begin{equation}\label{bee} \mathcal K \coloneq \Fun_k^{\mrm{L}}(\QC(X), \QC(X)) \times_{\Fun_k^{\mrm{L}}(\QC(X), \QC(X^{\mathscr E[1]}))} \Fun_k^{\mrm{L}}(\QC(X), \QC(X)) \end{equation}
may be identified with a triple $(\mrm{id}, \mrm{id}; \beth(\sigma))$ where $\beth(\sigma): u^* \simeq u^*$ is an equivalence of $k$-linear functors.

\begin{lemma}\label{pullback-lemma} Let $(\sigma; \alpha) \in \mrm{Map}_k'(\QC(X), \QC(X^\mathscr E))$ be a ($k$-linear) section $\sigma$ of $i^*$ with a specified equivalence $\alpha \colon i^* \sigma \xrightarrow{\simeq} \mrm{id}$.
There is a pullback square \[\begin{tikzcd}
	{(s_{0})_* \sigma} & {\mrm{id}} \\
	{\mrm{id}} & {u_* u^*}
	\arrow[from=1-1, to=1-2]
	\arrow[from=1-1, to=2-1]
	\arrow["{\overline{\beth(\sigma)}}", from=1-2, to=2-2]
	\arrow["{\mathrm{unit}}", from=2-1, to=2-2]
	\arrow["\lrcorner"{anchor=center, pos=0.125, rotate=180}, draw=none, from=2-2, to=1-1]
\end{tikzcd}\]	
in $\Fun_k^{\mrm{L}}(\QC(X), \QC(X))$, where $\overline{\beth(\sigma)}$ is adjoint to $\beth(\sigma) \colon u^* \simeq u^*$, and the other maps are the composites $(s_0)_* \sigma \xrightarrow{\mrm{unit}} (s_0)_* i_* i^* \sigma \simeq \mrm{id}$ (via $\alpha$ and the canonical equivalence $(s_0)_* i_* \simeq \mrm{id}$).
\end{lemma}
\begin{proof}
Consider the diagram \[\begin{tikzcd}
	{\QC(X)} && {\QC(X^\mathscr E)} \\
	& {\QC(X) \times_{\QC(X^{\mathscr E[1]})}\QC(X)}
	\arrow["{{s_0^* }}", from=1-1, to=1-3]
	\arrow["\tau"', from=1-1, to=2-2]
	\arrow["\iota", from=1-3, to=2-2]
\end{tikzcd}\]
of $k$-linear $\infty$-categories, where we recall that $\iota$ is fully faithful (\cref{cartesian}). 

The colimit-preserving $k$-linear functors $\iota$ and $\tau$ have right adjoints $\iota^R$ and $\tau^R$, which are canonically $k$-linear by \cref{adjoints} since $\iota$ and $\tau$ preserve compact objects (by another application of \cref{cartesian}). Applying $\Fun_k^{\mrm{L}}(\QC(X),-)$, we obtain a triangle of left-adjoint functors with $\mathcal K$ as the bottom vertex.

The map $\alpha$ supplies an equivalence $\iota \sigma \simeq (\mrm{id}, \mrm{id}; \beth(\sigma)\colon u^* \simeq u^*)$, and since $\iota$ is fully faithful,
\[(s_0)_* \sigma \simeq (s_0)_* \iota^R \iota \sigma \simeq \tau^R \iota \sigma \simeq \tau^R(\mrm{id}, \mrm{id}; \beth(\sigma)).\]
This final term is simply $\mrm{id} \times_{u_* u^*; \overline{\beth(\sigma)}} \mrm{id}$, and so we obtain our desired equivalence. 
\end{proof}
\begin{proof}[Proof of \cref{comparison}]
We use the notation of \cref{classify-sections} and its proof. Taking $\sigma = s_{\delta}^*$ in \ref{pullback-lemma} and using the identification $u_* u^* = \mrm{id} \oplus (-) \otimes \mrm{L}_{X/k}[1]$, we obtain a diagram
\[\begin{tikzcd}
	{(s_0)_* s_{\delta}^*} & {\mrm{id}} \\
	{\mrm{id}} & {u_* u^*} \\
	0 & {(-) \otimes \mrm{L}_{X/k}[1]}
	\arrow[from=1-1, to=1-2]
	\arrow[from=1-1, to=2-1]
	\arrow["\overline{\beth(\sigma)}", from=1-2, to=2-2]
	\arrow["\eta", from=2-1, to=2-2]
	\arrow[from=2-1, to=3-1]
	\arrow["\lrcorner"{anchor=center, pos=0.125, rotate=180}, draw=none, from=2-2, to=1-1]
	\arrow[from=2-2, to=3-2]   
	\arrow[from=3-1, to=3-2]
	\arrow["\lrcorner"{anchor=center, pos=0.125, rotate=180}, draw=none, from=3-2, to=2-1]
\end{tikzcd}\]
in $\Fun_{k}^{\mrm{L}}(\QC(X), \QC(X))$, where both squares are cocartesian. The maps $(s_{0})_* s_{\delta}^* \to \mrm{id}$ are those appearing in the construction of the classical Atiyah transformation, so this diagram exhibits the composite \[\mrm{id} \xrightarrow{\overline{\beth(\sigma)}} u_* u^* \to (-) \otimes \mrm{L}_{X/k}[1]\] as the classical Atiyah transformation. Comparing the description at the end of the proof of \cref{classify-sections}, this is also none other than the derived Atiyah transformation of \cref{derivedatiyah}.
\end{proof}

\subsection{The Atiyah class and Hochschild homology} \label{relating}

We now wish to relate the Atiyah class to Hochschild homology and complete the proof of \cref{HKR-atiyah}.

\begin{definition} Let $X$ be a qcqs derived scheme over $k$. The \emph{derived loop space} $\mathscr LX$ is the self-intersection of the diagonal $X \times_{X \times_k X} X$; we denote the two projections by $q_1, q_2 \colon \mathscr LX \rightrightarrows X$. The \emph{Hochschild homology sheaf} on $X$ is the quasicoherent sheaf $\underline{\mrm{HH}}_\bullet = \underline{\mrm{HH}}_{X/k, \bullet} = (q_2)_* \mathcal O_{\mathscr L X}$.  
\end{definition}

\begin{remark} (Comparison with other definitions)
\begin{enumerate} \item Suppose $X$ is a qcqs derived scheme over $k$. The assignment $U \mapsto \mathscr L U$ forms a Zariski cosheaf on $X$ by \cite[Lemma 4.2.]{LoopSpacesConnections}, so for any open affine $U = \Spec R \subset X$, we have a canonical equivalence $\underline{\mrm{HH}}_{X/k, \bullet}(U) \xrightarrow{\simeq} \underline{\mrm{HH}}_{U/k,\bullet}(U) = R \otimes_{R \otimes_k R} R$ of $R$-modules. In particular, we see that the quasicoherent sheaf $\underline{\mrm{HH}}_\bullet$ agrees with what one obtains by sheafifying the assignment of Hochschild homology on affines. 
\item Moreover, qcqs base change provides a canonical equivalence $\Delta^* \Delta_* \mathcal O_X \xrightarrow{\simeq} (q_1)_* (q_2)^* \mathcal O_{X} = \underline{\mrm{HH}}_\bullet$.
\end{enumerate}
\end{remark}

As usual, we equip $\underline{\mrm{HH}}_{\bullet}$ with the HKR filtration $\mathrm{Fil}^\star_{\mrm{HKR}} \underline{\mrm{HH}}_\bullet$, with associated graded $\mrm{gr}^i_{\mrm{HKR}} \underline{\mrm{HH}}_\bullet \simeq \Lambda^i \mrm{L}_{X/k}[i]$, by sheafifying the HKR filtration on affine opens.

\begin{notation} Given an object $M$ of an $\infty$-category $\mathcal C$, write $M(i)$ for the filtered object with $\mrm{Fil}^j M(i) = M$ for $j \ge i$ and $0$ otherwise.
\end{notation}

\begin{construction} \label{ci-maps} There are filtered maps $c_i \colon \mrm{Fil}^\star_{\mrm{HKR}} \underline{\mrm{HH}}_\bullet \to \Lambda^i \mrm{L}_{X/k}[i](i)$ obtained by using the maps $\mrm{HH}(R/k) \to \Omega^i_{R/k}[i]$ called $\pi_i$ in \cite[\S 1.3]{loday} on polynomial $k$-algebras $R$, and subsequently left Kan extending to all animated $k$-algebras and using Zariski descent. By \cite[Proposition 1.3.15]{loday}, the induced map $\Lambda^i \mrm{L}_{X/k}[i] \simeq \mrm{gr}^i_{\mrm{HKR}}\underline{\mrm{HH}}_\bullet \longrightarrow \Lambda^i \mrm{L}_{X/k}[i]$ is multiplication by $i!$. We call $c_1$ the HKR map; it induces an equivalence $\mrm{gr}^1_{\mrm{HKR}} \underline{\mrm{HH}} \xrightarrow{\simeq} \mrm{L}_{X/k}[1]$ in any characteristic.
\end{construction}

\begin{remark} \label{degen} The arguments of \cite[Lemma 3.5]{AV} and \cite[Theorem 4.8]{Yeku} immediately generalise to show that the maps $c_i$ induce an equivalence $\underline{\HH}_\bullet \xrightarrow{\simeq} \bigoplus_{i \in \N} \Lambda^i \mrm{L}_{X/k}[i]$ in characteristic $0$ and an equivalence $\underline{\HH}_\bullet/\mrm{Fil}^p_{\mrm{HKR}} \underline{\HH}_\bullet \xrightarrow{\simeq} \bigoplus_{i < p} \Lambda^i \mrm{L}_{X/k}[i]$ in characteristic $p$. A canonical partial splitting of the HKR filtration also appears as \cite[Theorem 3.2.1]{mundinger}; we expect these splittings to coincide.
\end{remark}

The following is the main ingredient in the proof of \cref{HKR-atiyah}: 
\begin{proposition} \label{HKR-Atiyah-Lemma} Let $X$ be a qcqs derived $k$-scheme. Then the map $\Delta_* \mathcal O_X \to \Delta_* \mrm{L}_{X/k}[1]$ of Fourier--Mukai kernels corresponding to the Atiyah transformation is adjoint to the HKR map \[\Delta^* \Delta_* \mathcal O_X \simeq \underline{\mrm{HH}}_\bullet \xrightarrow{\mrm{HKR}} \mrm{L}_{X/k}[1].\] 
\end{proposition}

\begin{remark} Under stronger assumptions (e.g. for $X$ quasicompact, separated, and smooth over $k$), \cref{HKR-Atiyah-Lemma} follows from \cite[Proposition 4.4]{Căldăraru} (and its proof); whilst Căldăraru works in characteristic $0$, which allows for statements concerning higher stages of the HKR filtration, this assumption is not needed to obtain the statement above. We will take the affine case of Căldăraru's proof and extend it to qcqs derived $k$-schemes. Note that a variant for (possibly nonsmooth) complex analytic spaces also appears as \cite[Theorem 5.1.3]{BuchweitzFlennerDecomp}.
\end{remark}

\begin{proof}[Proof of \cref{HKR-Atiyah-Lemma}]
  
We consider a slight variant $\mathscr L' X \coloneq X \times_{X \times_k X} X^{\mrm{L}_{X/k}}$ of the loop scheme, where the map $\Delta' = (s_0, s_\delta) \colon X^{\mrm{L}_{X/k}} \to X \times_k X$ is as in \cref{classicalatiyah}. This fits into the following diagram, where both squares are cartesian:
\[\begin{tikzcd}
	{\mathscr LX} & {\mathscr L'X} & X \\
	X & {X^{\mrm{L}_{X/k}}} & {X \times_k X}
	\arrow[from=1-1, to=1-2]
	\arrow["{q_1}", from=1-1, to=2-1]
	\arrow["{q_2'}"', from=1-2, to=1-3]
	\arrow["{q_1'}", from=1-2, to=2-2]
	\arrow["\Delta", from=1-3, to=2-3]
	\arrow["{i}", from=2-1, to=2-2]
	\arrow["{\Delta'}", from=2-2, to=2-3]
      \end{tikzcd}.\]
We recall that the adjoint of the Atiyah transformation is the composite
\[\Delta^* \Delta_* \mathcal O_X \simeq \Delta^* (\Delta')_* i_* \mathcal O_X \to \Delta^* (\Delta')_*i_* \mrm{L}_{X/k}[1] \simeq \Delta^* \Delta_* \mrm{L}_{X/k}[1] \to \mrm{L}_{X/k}[1],\]
where the second map is obtained from the fibre sequence $i_* \mrm{L}_{X/k} \to \mrm{sqz}_{\mathcal O_X}(\mrm{L}_{X/k}) \to i_* \mathcal O_X$.
By qcqs base change, we identify this composite with a map
\[\underline{\mrm{HH}}_\bullet = (q_2)_* (q_1)^* \mathcal O_X \simeq (q_2')_*(q_1')^*(i_* \mathcal O_X) \to (q_2')_* (q_1')^* (i_* \mrm{L}_{X/k}[1]) \simeq (q_2)_* q_1^* \mrm{L}_{X/k}[1] \to \mrm{L}_{X/k}[1];\]
this can be written more concretely as the composite
\[\phi_X \colon \mathcal O_{\mathscr L X} \simeq \mathcal O_{\mathscr{L}'X} \otimes_{\mrm{sqz}_{\mathcal O_X}(\mrm{L}_{X/k})}  \mathcal O_X \to \mathcal O_{{\mathscr L}'X} \otimes_{\mrm{sqz}_{\mathcal O_X}(\mrm{L}_{X/k})} \mrm{L}_{X/k}[1] \simeq \mathcal O_{\mathscr LX} \otimes_{\mathcal O_X} \mrm{L}_{X/k}[1] \to \mrm{L}_{X/k}[1],\] 
where we view $\mathcal O_{\mathscr LX}$ and $\mathcal O_{\mathscr L' X}$ as $\mathcal O_X$-bimodules via the projections to $X$, and the last map is induced by the inclusion $X \to \mathscr LX$.

The formation of $\phi_X$ is compatible with restrictions to open subschemes, since the assignments $U \mapsto \mathscr LU, \mathscr L'U$ form Zariski cosheaves. Evaluated on open affines $U = \Spec R \subseteq X$, this is the map of $R$-modules
\begin{align*} \phi_R \colon R \otimes_{R \otimes_k R} R \simeq &\,(R \otimes_{R \otimes_k R} \mrm{sqz}_R(\mrm{L}_{R/k})) \otimes_{\mrm{sqz}_R(\mrm{L}_{R/k})} R \\ \xrightarrow{\mrm{At}}&\,(R \otimes_{R \otimes_k R} \mrm{sqz}_R(\mrm{L}_{R/k})) \otimes_{\mrm{sqz}_R(\mrm{L}_{R/k})} \mrm{L}_{R/k}[1] \\ \simeq &\,(R \otimes_{R \otimes_k R} R) \otimes_{R} \mrm{L}_{R/k}[1] \\ \to& \mrm{L}_{R/k}[1],
\end{align*}
where the second map is again obtained from the boundary in the principal parts sequence. 
The formation of $\phi_R$ is functorial and, since this is true of its source and target, $R \mapsto \phi_R$ is the left Kan extension of its restriction to polynomial $k$-algebras. If $R$ is a smooth $k$-algebra, then the induced map $\mrm{HH}_\bullet(R/k) \to \Omega^1_{R/k}[1]$ is shown to be naturally homotopic to the HKR map in \cite[Proposition 4.4]{Căldăraru}. That this homotopy is natural in $R$ can be seen from the affine part of Căldăraru's proof, as a natural commutative diagram is exhibited on the chain level.
\end{proof}

 We can now finally prove the promised \cref{HKR-atiyah}:  
\begin{proof}[Proof of \cref{HKR-atiyah}.] The first equivalence of (1) is given by \cref{classify-sections}. For the second, we may combine the Fourier--Mukai equivalences recalled in \cref{fm} with adjunction to identify
\begin{align*}
\mrm{Map}_{\Fun_k^{\mrm{L}}(\QC(X),\QC(X))}(\mrm{id}, - \otimes \mathscr E[1]) & \simeq \mrm{Map}_{\QC(X \times_k X)}(\Delta_* \mathcal O_X, \Delta_* \mathscr E[1]) \\
                                                                               &\simeq \mrm{Map}_{\QC(X)}(\Delta^* \Delta_* \mathcal O_X, \mathscr E[1]) \\ &\simeq \mrm{Map}_{\QC(X)}(\underline{\HH}_{X/k, \bullet}, \mathscr E[1]).\end{align*} Taking $\mathscr E = \mrm{L}_{X/k}$ and using \cref{HKR-Atiyah-Lemma}, we obtain \hyperref[two]{(\ref*{two})}.
\end{proof}

\section{The period map}\label{period}
In this section, we finally construct the `commutative-to-noncommutative' period map \[\theta \colon \mrm{Def}_X \longrightarrow \mrm{CatDef}_{\QC(X)}\] attached to a derived $k$-scheme $X$.

\subsection{Construction}
The functor $\QC\colon \mrm{DSch}^{\mrm{op}} \to \mrm{Alg}(\PrL)$ induces a commutative diagram of $\infty$-categories 
\[\begin{tikzcd}
	{\mrm{CAlg}^{\mrm{an}}} && {\mrm{DSch}^{\op}} && {\Fun(\Delta^1, \mrm{DSch}^{\op})} \\
	{\mrm{CAlg}^{\an}} && {\mrm{Alg}(\PrL)} && {\Fun(\Delta^1, \mrm{Alg}(\PrL)).}
	\arrow["{\mrm{Spec}}", from=1-1, to=1-3]
	\arrow[equals, from=1-1, to=2-1]
	\arrow[from=1-3, to=2-3]
	\arrow["{{{\mrm{ev}_0}}}"', from=1-5, to=1-3]
	\arrow[from=1-5, to=2-5]
	\arrow["{\Mod_{(-)}}", from=2-1, to=2-3]
	\arrow["{{{\mrm{ev}_0}}}"', from=2-5, to=2-3]
\end{tikzcd}\]
There is a functor $\Fun(\Delta^1,\mrm{Alg}(\PrL)) \to \Mod(\PrL)$, compatible with the projections to $\mrm{Alg}(\PrL)$, which on objects takes a morphism $\mathscr C \to \mathscr D$ in $\mrm{Alg}(\PrL)$ to the presentably $\mathscr C$-linear $\infty$-category $\mathscr D$ \cite[Proposition 3.4.1.3]{HA}. We can therefore modify the above to a diagram \[\begin{tikzcd}
	{\mrm{CAlg}^{\mrm{an}}} && {\mrm{DSch}^{\op}} && {\Fun(\Delta^1, \mrm{DSch}^{\op})} \\
	{\mrm{CAlg}^{\an}} && {\mrm{Alg}(\PrL)} && {\Mod(\PrL).}
	\arrow["{\mrm{Spec}}", from=1-1, to=1-3]
	\arrow[equals, from=1-1, to=2-1]
	\arrow[from=1-3, to=2-3]
	\arrow["{{{\mrm{ev}_0}}}"', from=1-5, to=1-3]
	\arrow[from=1-5, to=2-5]
	\arrow["{\Mod_{(-)}}", from=2-1, to=2-3]
	\arrow[from=2-5, to=2-3]
\end{tikzcd}\]
Taking fibre products furnishes a triangle
\[\begin{tikzcd}
	\mrm{CAlg}^{\mrm{an}}\mrm{DSch}^{\mrm{op}} && \mrm{LinCat}^{\mrm{an}} \\
	& {\mrm{CAlg}^{\mrm{an}}}
	\arrow["g", from=1-1, to=1-3]
	\arrow["p"', from=1-1, to=2-2]
	\arrow["q^{\mrm{an}}", from=1-3, to=2-2]
\end{tikzcd}\]
where $p$ and $q^{\mrm{an}}$ are respectively the cocartesian fibrations appearing in \cref{scheme-construction} and \cref{an-cat}. 

By construction, $g$ sends a pair $(A, f \colon Y \to \Spec A)$ to the $\infty$-category $\QC(Y)$ equipped with the $A$-linear structure arising from $f^*$.

\begin{lemma}\label{squirrel} The functor $g$ sends $p$-cocartesian edges to $q^{\mrm{an}}$-cocartesian ones.
\end{lemma}
\begin{proof}
In light of the descriptions of cocartesian edges in \cref{cocart1} and \cref{cocart2}, the lemma amounts to the following:
\begin{itemize}
\item[($\star$)] For any cospan $\Spec B~\longrightarrow~\Spec A~\longleftarrow~Y$ of derived schemes, the $\QC(\Spec A)$-linear map $\QC(Y)~\longrightarrow~\QC(Y \times_A B)$ induces an equivalence \[\QC(Y) \otimes_{\QC(\Spec A)} \QC(\Spec B) \xlongrightarrow{\simeq} \QC(Y \times_{A} B).\] 
\end{itemize}
This follows from \hyperref[tensor]{Proposition \ref*{qcqs-cats} (\ref*{tensor})}.
\end{proof}
By \cref{squirrel}, we may (co)restrict $g$ to a map between the wide subcategories on cocartesian edges and take slices, as in the constructions of the deformation functors, to produce a map of left fibrations 
\[\begin{tikzcd}
	{\mrm{Defor}[X]} && {\mrm{Defor}[\QC(X)]} \\
	& {\mrm{CAlg}^{\mrm{an}}_{/k}}
	\arrow[from=1-1, to=1-3]
	\arrow[from=1-1, to=2-2]
	\arrow[from=1-3, to=2-2]
\end{tikzcd}\]
Straightening, this corresponds to a map $\theta^+\colon \mrm{Def}_X^{+} \longrightarrow \mrm{CatDef}^{+}_{\QC(X)}$.

\begin{definition} Let $X$ be a derived $k$-scheme. We may restrict $\theta^+$ to $\mrm{CAlg}^{\an, \mrm{art}}_{/k}$ to obtain a map \[\theta\colon \mrm{Def}_{X} \longrightarrow \mrm{CatDef}_{\QC(X)}\] of proximate formal moduli problems. We call this the (commutative-to-noncommutative) \emph{period map}. 
\end{definition}

\begin{variant} Restricting to $\mrm{CAlg}^{\mrm{an}, \mrm{art}}_{\Lambda \sslash k}$, we can consider the composite map $\mrm{Def}_{X,\Lambda} \to \mrm{CatDef}_{\QC(X), \Lambda} \to \mrm{CatDef}_{\QC(X),\Lambda}^{\land}$ to the associated formal moduli problem (see \cref{completion}), which will be denoted by $\theta^{\land}$.
\end{variant}

\begin{remark} If $X$ is qcqs, then $\QC(Y)$ is compactly generated for any infinitesimal deformation $Y$ of $X$ by combining \cref{qcqs-propagation} and \hyperref[qcqs-tcg]{Proposition \ref*{qcqs-cats} (\ref*{qcqs-tcg})}. It follows that the period map $\theta$ factors through the summand $\mrm{CatDef}_{\QC(X),\Lambda}^\omega \subset \mrm{CatDef}_{\QC(X), \Lambda}$ in this case.
\end{remark}

\subsection{Identifying the period map on tangent fibres}
By work of Brantner--Mathew \cite[Theorem 1.23]{BM}, there is an equivalence of $\infty$-categories \[\mrm{Moduli}_{\Lambda}^{\mrm{an}} \xlongrightarrow{\simeq} \mrm{Alg}_{\mrm{Lie}^{\pi}_{\Lambda, \Delta}}\] which sends a formal moduli problem $F$ to its tangent fibre $\mrm{T}_F$. The (completed) period map $\theta^{\land} \colon \mrm{Def}_{X,\Lambda} \to \mrm{CatDef}_{\QC(X),\Lambda}^{\land}$ constructed above therefore corresponds to a morphism $\mrm{T}_{\theta}$ of derived $(\Lambda,k)$-partition Lie algebras. We are now ready to prove our main result, concerning the underlying map of $k$-modules.

\begin{proposition}\label{identification}
Let $X$ be a qcqs derived scheme over a field $k$ and $\Lambda$ be a complete local Noetherian ring with residue field $k$. The (completed) commutative-to-noncommutative period map \[\theta^{\land}\colon \mrm{Def}_{X, \Lambda} \longrightarrow \mrm{CatDef}_{\QC(X), \Lambda}^{\land}\] induces a map $\mrm{T}_{\theta}\colon \mrm{R}\Gamma(X, \mathbf T_{X/k})[1] \to \HH^\bullet(X/k)[2]$ on tangent fibres which is equivalent to the dual {HKR} map. \end{proposition}

Here $\mathbf T_{X/k} = \mrm{map}_{\QC(X)}(\mrm{L}_{X/k}, \mathcal O_X)$ is the $\mathcal O_X$-dual of the cotangent complex.
\begin{proof}[Proof of \cref{identification}]
Let $V \in \Mod_{k, \ge 0}$. Taking loop spaces at the trivial deformations and using the identifications of \hyperref[looped-def]{\S \ref*{looped-def}} and \hyperref[loops-catdef]{\S \ref*{loops-catdef}}, we obtain a canonical square, functorial in $V \in \Mod_{k, \ge 0}$:
\[ \begin{tikzcd}
	{\Omega\mrm{Def}_X(\mrm{sqz}_k(V))} & {\mrm{Map}_{X\sslash k}(X^{V}, X)} \\
	{\Omega \mrm{CatDef}_{\QC(X)}(\mrm{sqz}_k(V))} & {\mrm{Map}_{k}'(\QC(X), \QC(X^V))}
	\arrow[ from=1-1, to=2-1]
	\arrow["\simeq",from=1-1, to=1-2]
	\arrow["\phi_V'", from =1-2, to=2-2]
	\arrow["\simeq", from=2-1, to=2-2]
\end{tikzcd} \]
where $\phi_V'$ is induced by applying $\QC(-)$. This map $\phi'_V$ is the $\mathscr E = \mathcal O_X \otimes_k V$-component of the natural transformation \[\phi_{\mathscr E} \colon \mrm{Map}_{X \sslash k}(X^{\mathscr E}, X) \to \mrm{Map}_k'(\QC(X), \QC(X^{\mathscr E}))\] of functors $\QC(X)_{\ge 0} \to \mathcal S$ induced by applying $\QC(-)$, and, using the identifications furnished by \cref{HKR-atiyah} and \hyperref[loops-catdef]{Proposition \ref*{loops-catdef}}, it sits in a square
\[\begin{tikzcd}
	{\mrm{Map}_{X\sslash k}(X^{\mathscr E}, X)} & {\mrm{Map}_{\QC(X)}(\mrm{L}_{X/k}, \mathscr E)} \\
	{\mrm{Map}_k'(\QC(X), \QC(X^{\mathscr E}))} & {\mrm{Map}_{\QC(X)}(\underline{\HH}_\bullet[-1], \mathscr E)}
	\arrow["\simeq", from=1-1, to=1-2]
	\arrow["\phi_{\mathscr E}"', from=1-1, to=2-1]
	\arrow[from=1-2, to=2-2]
	\arrow["\simeq", from=2-1, to=2-2]
\end{tikzcd}\]
natural in $\mathscr E \in \QC(X)_{\ge 0}$. By the Yoneda lemma, the right vertical map is given by composition with the map $h \colon \underline{\mrm{HH}}_{\bullet}[-1] \to \mrm{L}_{X/k}$ which is the image of $\mrm{id}_{\mrm{L}_{X/k}}$; by \cref{HKR-atiyah}, this is the HKR map.

We may now restrict and paste to obtain a diagram 
\[\begin{tikzcd}[column sep = 1.8em]
	{\mrm{Def}_{X}(\mrm{sqz}_k(V))} & {\Omega \mrm{Def}_X(\mrm{sqz}_k(V[1]))} & {\mrm{Map}_{\QC(X)}(\mrm{L}_{X/k}, \mathcal O_X \otimes_kV[1])} \\
	{\mrm{CatDef}_{\QC(X)}^{\land}(\mrm{sqz}_k(V))} & {\Omega \mrm{CatDef}_{\QC(X)}(\mrm{sqz}_k(V[1]))} & {\mrm{Map}_{\QC(X)}(\underline{\HH}_\bullet[-1], \mathcal O_X \otimes_k V[1])}
	\arrow["\simeq", from=1-1, to=1-2]
	\arrow[from=1-1, to=2-1]
	\arrow["\simeq", from=1-2, to=1-3]
	\arrow[from=1-2, to=2-2]
	\arrow["\mrm{HKR}", from=1-3, to=2-3]
	\arrow["\simeq", from=2-1, to=2-2]
	\arrow["\simeq", from=2-2, to=2-3]
\end{tikzcd}\]
functorial in $V \in \Mod_{k, \ge 0}^{\omega}$. Here we have identified 
\[\Omega \mrm{CatDef}_{\QC(X)}(A) = \Omega \mathrm{CatDef}_{\QC(X)}^\omega(A) = \Omega \mathrm{CatDef}_{\QC(X)}^{\land}(A),\]
where $A \coloneq \mrm{sqz}_k(V[1])$, as we may since $\mrm{CatDef}_{\QC(X)}^\omega$ is a $1$-proximate formal moduli problem (\cref{tame}) and $\mrm{CatDef}_{\QC(X)}^{\omega}(A) \subset \mrm{CatDef}_{\QC(X)}(A)$ is a union of components.

The result about tangent fibres now follows formally by using the natural equivalences \[\mrm{Map}_{\QC(X)}(\mathscr F, \mathscr G \otimes_k M) \simeq \mrm{Map}_{k}(M^\vee, \mrm{map}(\mathscr F, \mathscr G))\] for any $\mathscr F, \mathscr G \in \QC(X)$ and $M \in \Mod_k^{\omega}$, and the (defining) equivalences \[\mrm{Map}_{k}(V^{\vee}[-n], \mrm{T}_F) \simeq F(\mrm{sqz}_k(V[n]))\] for any formal moduli problem $F$.
\end{proof}

\subsection{Injectivity estimates} To make use of this period map $\theta$, it is important to have conditions under which the `derivative' $\mrm{T}_{\theta}$, or equivalently $\mrm{HKR}^\vee$, is injective on homotopy groups. In particular, injectivity on $\pi_{-1}$ implies injectivity on obstruction classes for (classical) deformations, whilst injectivity in more negative degrees is related to derived deformations and the corresponding partition Lie algebra being abelian. See \hyperref[apps]{\S \ref*{apps}} for more details.

The key injectivity estimates on the dual HKR map are as follows. Below, we will always write $\mrm{HKR}^\vee$ for the map $\mrm{R}\Gamma(X, {\mathbf T}_X)[-1] \to \mrm{HH}^\bullet(X/k)$; the map $\mrm{T}_{\theta}$ acquires a shift.

\begin{proposition}\label{inj-conditions} 
Let $X$ be a qcqs derived $k$-scheme.
\begin{enumerate} 
\item \label{general-estimate} The map $\mrm{HKR}^\vee \colon \mrm{R}\Gamma(X, {\mathbf T}_X)[-1] \to \mrm{HH}^\bullet(X/k)$ is injective on $\pi_i$ for \emph{all} $i$ if $\mrm{char}(k)=0$, and for $i \ge n-p$ if $\mrm{char}(k)=p > 0$ and $X$ is $n$-truncated. 
\item \label{strong-hkr} $\mrm{HKR}^\vee$ is injective on \emph{all} homotopy groups if the strong HKR theorem holds for $X$. 
\item \label{weak-CY} If $X$ is a weak Calabi--Yau variety and $\mrm{char}(k) = p > 0$, then $\mrm{HKR}^\vee$ is injective on $\pi_i$ if and only if all differentials out of $\{\mrm{E}^{1,1+d+i}_{r}\}_{r \ge 2}$ in the HKR spectral sequence vanish.
\end{enumerate}
\end{proposition}
\begin{proof}
	
  We first prove (\ref*{general-estimate}). We may assume $\mrm{char}(k) = p > 0$, as if $\mrm{char}(k) = 0$, then the strong HKR theorem holds and we may pass to (2). The HKR map factors as $\underline{\HH}_\bullet~\to~\underline{\mrm{HH}}_\bullet/\mrm{Fil}_{\mrm{HKR}}^{p} \underline{\HH}_\bullet \to \mrm{L}_{X/k}[1]$. The map on the right is canonically split (\cref{degen}), so it suffices to show that the natural map \[\phi \colon \mrm{map}_{\QC(X)}(\underline{\mrm{HH}}_\bullet/\mrm{Fil}_{\mrm{HKR}}^{p}\underline{\HH}_\bullet, \mathcal O_X) \longrightarrow \mrm{map}_{\QC(X)}(\underline{\HH}_\bullet, \mathcal O_X)\] of mapping spectra is injective on $\pi_i$ for $i \ge n-p$. Since $\mrm{Fil}^p_{\mrm{HKR}} \underline{\mrm{HH}}_\bullet$ is $p$-connective and $\mathcal O_X$ is $n$-truncated,
 $\mrm{fib}(\phi) = \mrm{map}_{\QC(X)}(\mrm{Fil}^{p}_{\mrm{HKR}} \underline{\HH}_\bullet[1], \mathcal O_X)$ is $(n-p-1)$-truncated and the result follows.

For (\ref*{strong-hkr}), we observe that if the strong HKR theorem holds for $X$, then in particular the HKR map admits a right inverse, so its dual has a left inverse. 

To prove (\ref*{weak-CY}), suppose that $X/k$ is a weak Calabi--Yau variety of dimension $d$. Having fixed a trivialisation of the canonical bundle, Serre duality provides a commutative square 
\[\begin{tikzcd}
	{\pi_i \mathrm{map}(\Omega^1_{X/k}[1], \mathcal O_X)} & {\pi_i\mathrm{map}(\underline{\mathrm{HH}}_\bullet, \mathcal O_X)} \\
	{\pi_{-d-i}\mathrm{R}\Gamma(X, \Omega^1_{X/k}[1])^\vee} & {\mathrm{{HH}}_{-d-i}(X/k)^\vee}
	\arrow["{\pi_i \mathrm{HKR}^\vee}", from=1-1, to=1-2]
	\arrow["\simeq"', from=1-1, to=2-1]
	\arrow["\simeq", from=1-2, to=2-2]
	\arrow[from=2-1, to=2-2]
\end{tikzcd}\]
where $(-)^\vee$ denotes the $k$-linear dual. It is hence equivalent to show that the HKR map $\HH_{-i-d} (X/k) \to \mrm{H}^{d+i}(X, \Omega^1_{X/k}[1])$ is surjective. Recall now that the HKR spectral sequence, which is simply the hypercohomology spectral sequence in this setting, has signature $\mrm{E}_2^{i,j} = \mrm{H}^{j}(X, \Omega^i_{X/k}) \Rightarrow \mathrm{HH}_{i-j}(X/k).$ On homotopy sheaves, $\mrm{HKR} \colon \underline{\mrm{HH}}_\bullet \to \Omega^1_{X/k}[1]$ is an equivalence on $\underline{\pi}_1$ (and otherwise $0$), so the induced map on the $\mrm{E}_2$ pages of the hypercohomology spectral sequences is simply projection onto the first row $\mrm{E}_2^{1,*} = \mrm{H}^*(X,\Omega^1_{X/k}[1])$. Surjectivity of $\mrm{HKR}_{-d-i}$ is therefore equivalent to asking that all differentials out of the bidegree $(1, d + i + 1)$ terms vanish, as claimed.
\end{proof}

\begin{proof}[Proof of {\hyperref[thm-A]{Theorem \ref*{thm-A}}}]

Combine \cref{identification} and \cref{inj-conditions}.
\end{proof}
\begin{remark}\label{petrov-remark}
In the proof of \cref{inj-conditions} (2), we only actually used the fact that the first piece of the HKR filtration splits off, but this \emph{prima facie} weaker hypothesis is in fact equivalent to the strong HKR theorem holding for $X$, and even $\underline{\mrm{HH}}_\bullet$ being free as a quasicoherent sheaf of animated $\mathcal O_X$-algebras; we thank A. Petrov for sharing the following argument with us. Any map $\alpha \colon \mrm{L}_{X/k}[1] \to \underline{\HH}_{\bullet}$ of quasicoherent sheaves extends to a map
\[\widetilde{\alpha} \colon \mrm{LSym}_{\mathcal O_X}^*(\mrm{L}_{X/k}[1]) \to \underline{\HH}_\bullet\] of animated quasicoherent $\mathcal O_X$-algebras. If $\alpha$ induces an isomorphism on $\underline{\pi}_1$, then since $\underline{\pi}_* \underline{\HH}_{\bullet} = \mrm{LSym}_{\mathcal O_X}^*(\underline{\pi}_1 \underline{\HH}_\bullet)$ by the affine HKR theorem, we find that $\widetilde{\alpha}$ is an equivalence on all homotopy sheaves, and hence is itself an equivalence
\end{remark}
\subsection{A criterion for HKR degeneration}
To apply \cref{inj-conditions}, it is desirable to know when the HKR spectral sequence degenerates in characteristic $p > 0$; this has been studied by various authors (see the introduction). Here we observe a cohomological condition for HKR degeneration when $k$ is perfect. A key input is the `commutative square' \cref{diagram} of spectral sequences associated to a smooth and proper $k$-scheme $X$, due to Antieau--Bragg \cite{AntieauBragg}, which makes use of the motivic filtration of \cite{BMSII}.
\begin{figure} 
  \centering
  \begin{tikzcd}
    &\HH_*(X/k)\arrow[Rightarrow]{dr}{\text{Tate}}&\\
    \mrm{H}^*(X,\Omega^*_X)\arrow[Rightarrow]{ur}{\text{HKR}}\arrow[Rightarrow]{dr}[swap]{\text{Hodge--de Rham}}&&\mrm{HP}_*(X/k)\\
    &\mrm{H}^*_{\mrm{dR}}(X/k)\arrow[Rightarrow]{ur}[swap]{\text{de Rham--HP}}
  \end{tikzcd}
  \caption{The `Hodge quartet'; borrowed from \cite[Figure 1]{AntieauBragg}.} \label{diagram}
\end{figure}
For dimension reasons, the HKR and Tate spectral sequences simultaneously degenerate if and only if the Hodge--de Rham and de Rham--HP spectral sequences do (see \cite[Remark 2.5]{AV} or \cite[Remark 3.6]{ABM}). This, along with a standard argument with Adams operations (inspired in this case by \cite{Elmanto}), leads to the following:

\hkrcriterion*
\begin{proof}
By the above discussion, it suffices to show that the de Rham--HP spectral sequence degenerates. Consider the corresponding crystalline--TP spectral sequence \[\mrm{E}_2^{s,t} = \mrm{H}^{s-t}_{\mrm{cris}}(X/W) \implies \mrm{TP}_{-s-t}(X)\] arising from the motivic filtration on $\mrm{TP}(X)$ constructed in \cite{BMSII}. By \emph{lo{c}.\ {c}it.} \S 9.4, $\mrm{TP}(X)$ carries Adams operations $\psi^\ell$ for each $\ell \in \mathbf Z_p^\times$ arising from the action of the $p$-completed circle on $\mrm{THH}^{\land}_p$, and $\psi^\ell$ acts with weight $\ell^n$ on $\mrm{gr}^n \mrm{TP}$; note that since $X$ is smooth and proper over $k$, $\mrm{THH}(X)$ and $\mrm{TP}(X)$ are automatically $p$-complete. As the operations intertwine the differentials and the crystalline cohomology is $W$-free, all differentials in the crystalline--TP spectral sequence vanish; indeed, for each $\ell$, each $y = d_r(x)$ satisfies an equation of the form $\ell^i y = \ell^j y$ for $i \ne j$, and is hence torsion by picking $\ell$ not to be a root of unity. We hence obtain a noncanonical splitting $\mrm{TP}_{n}(X) \simeq \bigoplus_{i \in \mathbf Z} \mrm{H}^{n+2i}_{\mrm{cris}}(X/W)$.
Now $\mrm{HP}(X/k) \simeq \mrm{TP}(X)/p$ by (for example) \cite[Theorem 3.4]{Antieau-Mathew-Nikolaus}, and $\mrm{R}\Gamma_{\mrm{cris}}(X/W)/p \simeq \mrm{R}\Gamma_{\mrm{dR}}(X/k)$. Reducing the above splitting mod $p$ (and again using the torsion-free hypothesis), we obtain an analogous splitting $\mrm{HP}_n(X) = \bigoplus_{i \in \mathbf Z} \mrm{H}^{n+2i}_{\mrm{dR}}(X/k)$, and the de Rham--HP spectral sequence degenerates for dimension reasons.
\end{proof}
\section{Applications}\label{apps}

This explicit description of the period map and the study of its injectivity on homotopy groups have immediate deformation-theoretic consequences, which we now explore.

\subsection{Square-zero extensions and derived invariance}\label{lieblich-section}

We first provide a general statement about lifting along square-zero extensions of animated rings.
\begin{proposition}\label{lifts} Let $X$ be a qcqs derived scheme over a field $k$. Let $\mathfrak X/B$ be a lift of $X/k$ over some $B \in \mrm{CAlg}^{\mrm{an},\mrm{art}}_{\Lambda \sslash k}$, defining a corresponding lift $\QC(\mathfrak X)$ of $\QC(X)$.
	
Suppose $A \twoheadrightarrow B$ is a square-zero extension by $V[i-1]$ for some $V \in \Mod_k^{\heartsuit}$ and $i \ge 1$. If $\mrm{T}_{\theta}$ is injective on $\pi_{-i}$, then $\mathfrak X$ lifts along $A \to B$ if and only if $\QC(\mathfrak X)$ does.  
\end{proposition}
\begin{remark} Slight care is needed since $\mrm{CatDef}_{\QC(X)}$ is not necessarily a formal moduli problem; otherwise, this result is standard. \end{remark}
\begin{proof}[Proof of \cref*{lifts}]
The `only if' direction being clear, we consider the `if' direction. In general, if $F$ is a formal moduli problem, then the obstruction to lifting a point $x \in \pi_0 F(B)$ along $A \to B$ is the image of $x$ in \[\pi_0 F(\mrm{sqz}_k(V[i])) \simeq \pi_{-i}(\mrm{T}_F) \otimes_k V , \] since $F(A) \simeq F(B) \times_{F(\mrm{sqz}_k(V[i]))} F(k)$ (see e.g.\ \cite[\S 1.6.1.]{pridham} or \cite[\S 3.3]{BT}). Taking $F \coloneq \mrm{Def}_X$, we wish to show the obstruction to lifting $[\mathfrak X] \in \pi_0 \mrm{Def}_X(B)$ vanishes; this can be checked after applying $\mrm{T}_{\theta}$ since this map is injective on $\pi_{-i}$. Hence it is enough to show that the image $\theta^{\land}([\mathfrak X]) \in \pi_0 \mrm{CatDef}^{\land}_{\QC(X)}(B)$ of $[\mathfrak X]$ lifts. But $\theta^{\land}$ factors through $\mrm{CatDef}_{\QC(X)}$, and $\theta([\mathfrak X])$ lifts to $A$, so its image in $\mrm{CatDef}_{\QC(X)}^{\land}(A)$ is a lift of $\theta^{\land}([\mathfrak X])$.
\end{proof}

This almost immediately implies \hyperref[thm-B]{Theorem \ref*{thm-B} (1)} by restricting to discrete rings; we defer the proof to the end of this section. We can now prove the corollary:
\begin{proof}[Proof of \cref{def-iso}]

We first prove \hyperref[cor1]{(\ref*{cor1})}. It suffices to prove that $\pi_i \mrm{T}_{\theta}$ is injective for $i = -1$, bijective when $i = 0$, and surjective when $i = 1$; see \cite[Theorem 3.1.]{Manetti-Exposition} or \cite[Theorem 1.45]{pridham}.

We have already shown that $\pi_{-1} \mrm{T}_{\theta}$ is injective when $\mrm{char}(k) \ne 2$. When $\mrm{char}(k)=2$, we can upgrade the conclusion of \hyperref[general-estimate]{Theorem \ref*{thm-A} (\ref*{general-estimate})}, as at the last step of the proof, we have $\pi_{-3} \fib(\phi) \simeq \pi_0 \mrm{map}_{\QC(Z)}(\mrm{Fil}^2_{\mrm{HKR}} \underline{\mrm{HH}}_\bullet[-2], \mathcal O_Z) \simeq \pi_0 \mrm{map}_{\QC(Z)}(\Omega^2_{Z/k}, \mathcal O_Z) = \mrm{H}^0(Z, \Lambda^2 \mrm{T}_{Z/k}) = 0$, where the second equivalence uses the fact that $\mrm{Fil}^3_{\mrm{HKR}}[-2]$ is $1$-connective whilst $\mathcal O_Z$ is $0$-truncated.

The map $\pi_0 \mrm{T}_{\theta} = \mrm{HKR}^\vee \colon \mrm{H}^1(Z, \mrm{T}_{Z/k}) \to \mrm{HH}^2(Z/k)$ on tangent spaces is always injective by \hyperref[thm-A]{Theorem \ref*{thm-A}}, and it is surjective in this case by the vanishing hypotheses; the other contributions to $\mrm{HH}^2$ from the (dual) HKR spectral sequence $\mrm{H}^{i}(X, \Lambda^j \mrm{T}_{X/k}) \Rightarrow \mrm{HH}^{i+j}(X/k)$ vanish. The map $\pi_1 \mrm{T}_{\theta}$ is an isomorphism for the same reason; this completes the proof of \hyperref[cor1]{(\ref*{cor1})}.

The first part of \hyperref[cor2]{(\ref*{cor2})} follows immediately from \hyperref[cor1]{(\ref*{cor1})}. Finally, if $Z \in \{X,Y\}$ is projective, then the hypothesis $\mathrm{H}^{2}(Z, \mathcal O_Z) = 0$ means there is no obstruction to lifting an ample line bundle on $Z$ to any formal lift $\mathfrak Z/\mathrm{Spf}(\Lambda)$ (e.g. \cite[Lemma 08VR]{stacks-project}). Hence all formal lifts are uniquely algebraisable by the Grothendieck existence theorem \cite[Lemma 0A4Z, Theorem 089A]{stacks-project} and we are done.
\qedhere

\end{proof}

\begin{remark}
As mentioned in the introduction, \hyperref[thm-B]{Theorem \ref*{thm-B}} (1) and its corollary imply that liftability along square-zero extensions is a derived invariant (with some additional conditions when $p=2$). We are unsure as to whether this positive result should generalise to higher-order deformations: Addington--Bragg have found derived-equivalent Calabi--Yau threefolds $X$ and $M$ over $\overline{\mathbf F}_3$ where the deformation spaces have different dimensions $h^1(T_X) \ne h^1(T_M)$ \cite[Theorem 1.1.]{Add-Bragg}, so it seems reasonable that derived-equivalent Calabi--Yau varieties might display quite different behaviour as far as liftability is concerned.
\end{remark}

\begin{remark} (Infinitesimal categorical Torelli). Following \cite{toenAnel}, injectivity on $\pi_0$ of the period map  can be viewed as a kind of Torelli theorem. To\"{e}n--Anel prove such a result for smooth and proper schemes over an algebraically closed field [\citelink{toenAnel}{Ibid}, Proposition 5.1] via different methods, but also propose that it should hold because of the map on tangent spaces respecting the HKR decomposition [\citelink{toenAnel}{Ibid}, \S 8.3]. 
\end{remark}
\begin{remark} \label{derEquiv} (Derived equivalences).  There are multiple notions of derived equivalence for $k$-schemes, but most coincide for qcqs schemes. Indeed, suppose $X$ and $Y$ are $k$-schemes. If $\QC(X) \simeq \QC(Y)$ as $k$-linear $\infty$-categories, then the homotopy categories $\mrm{D}_{\mrm{qc}}(-)$ are equivalent as $k$-linear triangulated categories (simply by passing to homotopy categories); the converse holds when $X$ and $Y$ are qcqs, as then these triangulated categories possess unique dg enhancements \cite[Theorem B]{CNS}. Note that in the smooth and proper case, any $k$-linear exact equivalence is induced by a Fourier--Mukai transform by Olander's generalisation of Orlov's theorem \cite[Theorem 2.2]{Orlov} \cite{olander} and hence lifts to an equivalence of $k$-linear $\infty$-categories.
\end{remark}

\subsection{Abelian partition Lie algebras and derived-unobstructedness}\label{der-unob}

Let us now connect the period map to the notion of \emph{derived}-unobstructedness and abelianity of partition Lie algebras. Since this notion does not quite appear in the literature (outside characteristic $0$), we briefly study this notion, at first only in the equicharacteristic case (i.e.\ when $\Lambda = k$). In the following, we follow the conventions of \cite{BM}. 

\begin{construction}\label{aug-constr} (Augmenting the partition Lie algebra monad)
  
Consider the following pair of adjunctions:
\[\begin{tikzcd}
	{\mathrm{Mod}_{k, \ge 0}^{\mathrm{op}}} && {(\mathrm{CAlg}_k^{\mathrm{an,aug}})^{\mathrm{op}}} && {\mathrm{Mod}_k.}
	\arrow[""{name=0, anchor=center, inner sep=0}, "{{\mathrm{free}}}", shift left=2, from=1-1, to=1-3]
	\arrow[""{name=1, anchor=center, inner sep=0}, "{{\mathrm{forget}}}", shift left=2, from=1-3, to=1-1]
	\arrow[""{name=2, anchor=center, inner sep=0}, "{{\mathrm{cot}^\vee}}", shift left=2, from=1-3, to=1-5]
	\arrow[""{name=3, anchor=center, inner sep=0}, "{{\mathrm{sqz} \circ  \tau_{\ge 0} \circ (-)^\vee}}", shift left=2, from=1-5, to=1-3]
	\arrow["\dashv"{anchor=center, rotate=90}, draw=none, from=1, to=0]
	\arrow["\dashv"{anchor=center, rotate=90}, draw=none, from=3, to=2]
\end{tikzcd}\]
The composite adjunction is $\tau_{\ge 0} \circ (-)^\vee \dashv (-)^\vee$, and the adjunction induces a map of monads on $\Mod_k$ \[\widetilde{T^\vee} \longrightarrow (-)^\vee \circ \tau_{\ge 0} \circ (-)^{\vee},\] where $\widetilde{T^\vee}$ sends $M \mapsto \mrm{cot}(\mrm{sqz}(\tau_{\ge 0}(M^\vee)))^\vee$ (cf.\ \cite[Remark 4.51]{BM}). We may restrict these monads to $\Mod_{k, \le 0}^{\mrm{ft}}$ and form their sifted-colimit-preserving extensions to $\Mod_k$ using \cite[Corollary 3.18]{BM}, which are respectively $T^\vee = \mrm{Lie}^{\pi}_{k, \Delta}$ (by definition) and $\mrm{id}$. Since these sifted-colimit-preserving extensions are formed by left Kan extension, there is an induced commutative diagram of monads on $\Mod_k$ as follows:
\begin{equation}\label{monad-square}
\begin{tikzcd}
	{\mathrm{Lie}^{\pi}_{k, \Delta}} & {\mathrm{id}} \\
	{\widetilde{T^\vee}} & {(-)^\vee \circ \tau_{\ge 0} \circ (-)^{\vee}.}
	\arrow[from=1-1, to=1-2]
	\arrow[from=1-1, to=2-1]
	\arrow[from=1-2, to=2-2]
	\arrow[from=2-1, to=2-2]
\end{tikzcd}
\end{equation}
\end{construction}

\begin{definition} (Abelian partition Lie algebras, I) 
\begin{enumerate} \item The augmentation $\mrm{Lie}^{\pi}_{k, \Delta} \to \mrm{id}$ produced in \cref*{aug-constr} induces a functor $\mrm{triv} \colon \Mod_k \to \mrm{Lie}^{\pi}_{k, \Delta}$ which we call the \emph{trivial derived partition Lie algebra functor}. Composing with the equivalence $\mrm{Lie}^{\pi}_{k, \Delta} \simeq \mrm{Moduli}_k^{\mrm{an}}$, we obtain a functor $\mrm{triv}^{\mrm{FMP}} \colon \Mod_k \to \mrm{Moduli}_k$.
\item A partition Lie algebra is called \emph{abelian} if it lies in the essential image of $\mrm{triv}$. A derived formal moduli problem over $k$ is called \emph{derived-unobstructed} if it lies in the essential image of $\mrm{triv}^{\mrm{FMP}}$ (i.e.\ if the corresponding partition Lie algebra is abelian).
\end{enumerate} 
\end{definition}

The functor $\mrm{triv}^{\mrm{FMP}}$ admits a relatively straightforward description:
\begin{proposition} 
For each $M \in \Mod_k$, there is an equivalence $\mrm{triv}^{\mrm{FMP}}(M) \simeq \Omega^\infty(M \otimes_k \mathfrak m_{(-)})$, where $\mathfrak m_A \coloneq \mrm{fib}(A \to k)$ denotes the augmentation ideal.
\end{proposition}
\begin{proof}
Let us denote the functor $M \mapsto \Omega^\infty(M \otimes_k \mathfrak m_{(-)})$ by $\overline{\mrm{triv}}^{\mrm{FMP}}$, and write $\Psi$ for the equivalence $\mrm{Lie}_{k,\Delta}^{\pi} \xrightarrow{\simeq} \mrm{Moduli}_{k}^{\mrm{an}}$ sending $\mathfrak g \mapsto \mrm{Map}_{\mrm{Alg}_{\mrm{Lie}^{\pi}_{k, \Delta}}}(\mathfrak D(-), \mathfrak g)$. We will show that the functors $\Psi \circ \mrm{triv},  \overline{\mrm{triv}}^{\mrm{FMP}} \colon \Mod_k \to \mrm{Moduli}_{k}^{\mrm{an}}$ are equivalent. 

Using the perspective on the Koszul duality functor $\mathfrak D$ in \cite[Remark 4.51]{BM}, the square of monads \hyperref[monad-square]{(\ref{monad-square})} induces a commutative diagram
\[\begin{tikzcd}
	{\mathrm{Mod}_{k, \ge 0}^{\mathrm{op}}} && {(\mathrm{CAlg}^{\mathrm{an,\mathrm{aug}}}_k)^{\mathrm{op}}} \\
	{\mathrm{Mod}_k} && {\mathrm{Alg}_{\mathrm{Lie}^{\pi}_{k, \Delta}}} \\
	& {\mathrm{Mod}_k.}
	\arrow["{\mathrm{free}}", from=1-1, to=1-3]
	\arrow["{(-)^\vee}"', from=1-1, to=2-1]
	\arrow["{\mathfrak D}", from=1-3, to=2-3]
	\arrow["{\mathrm{triv}}", from=2-1, to=2-3]
	\arrow["{\mathrm{id}}"', from=2-1, to=3-2]
	\arrow[from=2-3, to=3-2]
\end{tikzcd}\]
Recall that $\mathfrak D \colon (\mrm{CAlg}_k^{\mrm{an}, \mrm{aug}})^{\mrm{op}} \to \mrm{Lie}_{k,\Delta}^{\pi}$ has a left adjoint $C^*$ (the \emph{Chevalley--Eilenberg cochains}), and that the counit of this adjunction is an equivalence when restricted to complete local Noetherian animated $k$-algebras. We find equivalences, natural in $M \in \mathrm{Mod}_k$, $A \in \mathrm{CAlg}_{k \sslash k}^{\mathrm{an,art}}$:
\begin{align*}
\mrm{Map}_{\mathrm{Alg}_{\mrm{Lie}^{\pi}_{k, \Delta}}}(\mathfrak D(A), \mathrm{triv}(M^\vee)) &\simeq \mathrm{Map}_{\mathrm{Alg}_{\mathrm{Lie}^{\pi}_{k, \Delta}}}(\mathfrak D(A), \mathfrak D(\mathrm{free}(M))) 
\\ &\simeq\mathrm{Map}_{\mathrm{Mod}_k}(M, \mathrm{forget}(C^*(\mathfrak D(A)))) 
\\ &\simeq\mathrm{Map}_{\mathrm{Mod}_k}(M, \mathrm{forget}(A)) 
\\ &\simeq \Omega^\infty ( M^\vee \otimes_k \mathfrak m_A),
\end{align*}
where we have used that for Artinian $A$, the unit $A \to C^*(\mathfrak D(A))$ is an equivalence (see \cite[Proposition 4.52]{BM}).

In other words, we have an equivalence $\Psi \circ \mrm{triv} \circ (-)^\vee \simeq \overline{\mrm{triv}}^{\mrm{FMP}} \circ (-)^\vee$ on $\Mod_{k, \ge 0}$. Composing with linear duality, we find an equivalence $\Psi \circ \mrm{triv}|_{\Mod_{k, \le 0}^{\mrm{ft}}} \simeq \overline{\mrm{triv}}^{\mrm{FMP}}|_{\Mod_{k, \le 0}^{\mrm{ft}}}$. We now claim that each of these functors on $\Mod_{k}$ preserves sifted colimits; then the equivalence on $\Mod_{k, \le 0}^{\mrm{ft}}$ extends to all of $\Mod_k$ by \cite[Proposition 3.14]{BM}, and we are done. Since $\Psi$ is an equivalence and $\overline{\mrm{triv}}^{\mrm{FMP}}$ preserves sifted colimits, it suffices to show that $\mrm{triv} \colon \Mod_k \to \mrm{Alg}_{\mrm{Lie}^{\pi}_{k, \Delta}}$ does. But the forgetful functor $U \colon \mrm{Alg}_{\mrm{Lie}^{\pi}_{k, \Delta}} \to \Mod_k$ reflects sifted colimits (as $\mrm{Lie}^{\pi}_{k, \Delta}$ preserves them), and the composite $U \circ \mrm{triv}$ is the identity, so we are done.
\end{proof} 

Using this moduli-theoretic description, we can generalise the notions of derived-unobstructed formal moduli problems and abelian partition Lie algebras, and relax the assumption that $\Lambda = k$: 
\begin{definition}\label{unob-def} (Abelian partition Lie algebras, \textsc{II})

\begin{enumerate}
\item An $\mathbf E_0$-$(\Lambda,k)$ formal moduli problem is called \emph{derived-unobstructed} if it is equivalent to one of the form
\[ \mrm{triv}^{\mrm{FMP}}(M) \colon B \mapsto \Omega^{\infty}(M \otimes_\Lambda \mathfrak{m}_B)\] for some graded-free $\Lambda$-module $M$, where $\mathfrak m_B = \mrm{fib}(B \to k)$ denotes the augmentation ideal. 
\item More generally, an object of $\mrm{Moduli}_{\Lambda}^{\mrm{an}}$ or $\mrm{Moduli}_{\Lambda}^{\mathbf E_n}$ is called derived-unobstructed if it is the restriction of a derived-unobstructed $\mathbf E_0$-formal moduli problem.
\item A derived (or spectral) $(\Lambda,k)$-partition Lie algebra is called \emph{abelian} if the corresponding derived (resp. spectral) $(\Lambda,k)$-formal moduli problem is derived-unobstructed. 
We write $\mrm{triv}(M)$ for the partition Lie algebra corresponding to (the appropriate restriction of ) the formal moduli problem $\mrm{triv}^{\mrm{FMP}}(M)$.
\end{enumerate}
\end{definition}

\subsubsection{Some basic theory of abelian partition Lie algebras}
We now collect some straightforward facts about abelian partition Lie algebras:
\begin{proposition} Let $F \in \mrm{Moduli}_{\Lambda}^{\mrm{an}}$ be a derived-unobstructed derived formal moduli problem, and $M$ be a graded-free $\Lambda$-module.
\begin{enumerate}
	\item \label{ab-a} The underlying $k$-module of $\mrm{triv}(M)$ is $M \otimes_{\Lambda} k$.
	\item \label{ab-b} The underlying classical deformation functor $F^{\mrm{cl}} = \pi_0 F \colon \mathrm{CRing}^\mrm{art}_{\Lambda \sslash k} \longrightarrow \mathrm{Set}$ of $F$ is given by $B \mapsto \pi_0(\mrm{T}_F) \otimes_\Lambda \mathfrak m_B$, and is formally smooth/unobstructed in the usual sense of \emph{\cite[Definition 2.2]{Schlessinger}}.
	\item \label{ab-c} If $B \to C$ is any map in $\mathrm{CAlg}^{\mathrm{art}}_{\Lambda \sslash k}$ surjective on $\pi_*$, then $F(B) \to F(C)$ is surjective on $\pi_*$.
	\item \label{ab-d} If $M$ is a \emph{finite type} $\Lambda$-module, then there is a natural equivalence \[\mrm{triv}^{\mrm{FMP}}(M)~\simeq~\mrm{Map}_{\mathrm{DAlg}_{\Lambda \sslash k}}(\mrm{LSym}^*_\Lambda(M^\vee),-).\] 

\end{enumerate}
\end{proposition}

Here we write $\mathrm{DAlg}$ for the $\infty$-category of \emph{derived rings} as developed by Bhatt--Mathew, Brantner--Mathew, and Raksit (and ultimately going back at least to Illusie); see \cite[\S 4]{raksit-2020}. These can be modelled by simplicial-cosimplicial rings \cite[Corollary 5.29]{BCN}.

\begin{proof}
For (\ref*{ab-a}), recall that the underlying $k$-module of $\mrm{triv}(M)$ is determined by there being equivalences $\Omega^\infty(V[n] \otimes_k \mrm{triv}(M)) \simeq \mrm{triv}^{\mrm{FMP}}(M)(\mrm{sqz}_k(V[n]))$ (natural in $V \in \Mod_k^{\omega, \heartsuit}$). This second space is $\Omega^\infty(M \otimes_{\Lambda} V[n]) \simeq \Omega^\infty( (M \otimes_{\Lambda} k) \otimes_k V[n])$, and so we find that $\mrm{triv}(M) \simeq M \otimes_\Lambda k$.

The first part of (\ref*{ab-b}) is clear, and the second part is a special case of (\ref*{ab-c}). Statement (\ref*{ab-c}) holds by viewing $\mrm{triv}^{\mrm{FMP}}(M)$ as the composite of three functors, each of which preserves $\pi_*$-surjective maps. For (\ref*{ab-d}), we compute 
$\mrm{triv}^{\mrm{FMP}}(M)(B) \simeq \mathrm{Map}_\Lambda(M^\vee, \mathfrak m_B) 
\simeq \mathrm{Map}_{\Lambda \sslash k}(\mathrm{LSym}^*_\Lambda(M^\vee), B)$. 
\end{proof}

\begin{remark}
The notion of derived-unobstructedness here does \emph{not} imply (formal) smoothness in the sense of \cite[Proposition 6.2.8]{DAG}, \cite[Definition 1.27]{pridham}. Indeed, a formal moduli problem is smooth if and only if the tangent fibre is connective (\cite[Proposition 6.2.8]{DAG}, \cite[Remark 1.28]{pridham}). However, formal moduli problems smooth in this strong sense seem rare in practice.
\end{remark}

A standard way to show that a DGLA is homotopy abelian is to find a map to an abelian DGLA which is injective on homotopy groups; cf.\ \cite[Proposition 4.11 (ii)]{kkp}. This result generalises to partition Lie algebras by modifying an argument we learned from \cite[Lemma 2.7]{FM18}:

\begin{proposition} \label{reflect-abelian} Suppose $i \colon \mathfrak g \to \mathfrak h$ is a map of $(\Lambda,k)$-partition Lie algebras which is injective on homotopy groups. If $\mathfrak h$ is abelian, then so is $\mathfrak g$.
\end{proposition}

\begin{proof}
We may assume $\mathfrak h = \mrm{triv}(M)$ for some graded-free $\Lambda$-module $M$, and will write $\mrm{forget}$ for the forgetful functor $\mrm{Alg}_{\mrm{Lie}^{\pi}_{\Lambda, \Delta}} \to \Mod_{k}$. Since $\mrm{forget}(\mathfrak g) \to \mrm{forget}(\mathfrak h) = M \otimes_\Lambda k$ is injective on homotopy groups, and the modules are graded-free, we can find a left inverse $r \colon M \otimes_\Lambda k \to \mrm{forget}(\mathfrak g)$ to $\mrm{forget}(i)$. We may pick an equivalence $\mrm{forget}(\mathfrak g) \simeq N \otimes_\Lambda k$ for a graded-free $\Lambda$-module $N$, and lift $r$ arbitrarily to a map $M \to N$. Then the composite $\mathfrak g \xlongrightarrow{i} \mrm{triv}(M) \longrightarrow \mrm{triv}(N)$ is an equivalence, since it is so on underlying $k$-modules.
\end{proof}
Combining this with previous results, we now have a proof of \hyperref[thm-B]{Theorem \ref*{thm-B}}:

\begin{proof}[Proof of Theorem \ref*{thm-B}.]
Part (1) follows by taking $i=1$ in \cref{lifts}; $\mrm{T}_{\theta}$ is injective on $\pi_{-1}$ by \hyperref[thm-A]{Theorem \ref*{thm-A}} (2) when $\mrm{char}(k) \ne 2$, and otherwise by assumption. Part (\ref{abelianity}) is a special case of \cref{reflect-abelian}.   
\end{proof} 

\emergencystretch 1em
\renewcommand*{\bibfont}{\Small}
\printbibliography
\end{document}